\documentclass[12pt]{amsart}
\usepackage{amsmath,amsthm,amsfonts,amssymb,graphicx,psfrag,dsfont}

\usepackage{hyperref}
\usepackage[alphabetic,initials]{amsrefs}

\psfrag{tau0}{$\tau=0$}
\psfrag{tau1}{$\tau=\frac{1}{2} - \epsilon$}
\psfrag{tau2}{$\tau=\frac{1}{2}$}
\psfrag{tau3}{$\tau=\frac{1}{2} + \epsilon$}

\theoremstyle{plain}
\newtheorem*{thmu}{Theorem}

\newtheorem{thm}{Theorem}
\newtheorem{prop}{Proposition}[section]

\newtheorem{cor}[prop]{Corollary}
\newtheorem{lemma}[prop]{Lemma}

\theoremstyle{definition}

\newtheorem{defn}[prop]{Definition}
\newtheorem*{notation}{Notation}

\theoremstyle{remark}
\newtheorem{remark}[prop]{Remark}
\newtheorem*{remarku}{Remark}

\newcommand{\asympdef}{\operatorname{def}_\infty}

\newcommand{\const}{^{\operatorname{c}}}
\newcommand{\core}{^{\operatorname{K}}}
\newcommand{\cov}{{\operatorname{cov}_{\infty}}}

\newcommand{\HC}{\operatorname{HC}}
\newcommand{\Hom}{\operatorname{Hom}}

\newcommand{\ind}{\operatorname{ind}}

\newcommand{\interior}{\operatorname{int}}

\newcommand{\loc}{\operatorname{loc}}

\newcommand{\muCZ}{\mu_{\operatorname{CZ}}}
\newcommand{\piTu}{\pi Tu}

\newcommand{\trivial}{^{\operatorname{t}}}

\newcommand{\wind}{\operatorname{wind}}
\newcommand{\windpi}{\operatorname{wind}_\pi}

\newcommand{\CC}{{\mathbb C}}

\newcommand{\NN}{{\mathbb N}}

\newcommand{\RR}{{\mathbb R}}
\newcommand{\ZZ}{{\mathbb Z}}

\newcommand{\dD}{{\mathcal D}}

\newcommand{\fF}{{\mathcal F}}
\newcommand{\hH}{{\mathcal H}}

\newcommand{\tT}{{\mathcal T}}
\newcommand{\uU}{{\mathcal U}}

\newcommand{\p}{\partial}
\newcommand{\Cinftyloc}{C^\infty_{\loc}}
\newcommand{\oC}{\operatorname{C}}
\newcommand{\oE}{\operatorname{E}}
\newcommand{\oN}{\operatorname{N}}

\numberwithin{equation}{section}

\title[Compactness for holomorphic curves in $3$--manifolds]{Compactness for
Embedded Pseudoholomorphic Curves in $3$--manifolds}
\author{Chris Wendl}
\address{ETH Z\"urich \\
Departement Mathematik, HG G38.1 \\ 
R\"amistrasse 101 \\
8092 Z\"urich \\ 
Switzerland}
\email{wendl@math.ethz.ch}
\urladdr{http://www.math.ethz.ch/~wendl/}
\thanks{Research partially supported by an NSF Postdoctoral Fellowship
(DMS-0603500) and a DFG grant (CI 45/2-1)}

\begin{document}

\begin{abstract}
We prove a compactness theorem for holomorphic curves in $4$--dimensional 
symplectizations that have embedded projections to the 
underlying $3$--manifold.
It strengthens the cylindrical case of the SFT compactness theorem 
\cite{SFTcompactness} by using intersection theory to show that degenerations
of such sequences never give rise to multiple covers or nodes, so
transversality is easily achieved.  This has application to the theory
of stable finite energy foliations introduced in \cite{HWZ:foliations}, 
and also suggests a new approach to defining
SFT-type invariants for contact $3$--manifolds, or more generally, 
$3$--manifolds with stable Hamiltonian structures.
\end{abstract}

\maketitle

\tableofcontents

\section{Introduction and main results}
\label{sec:intro}

Compactness arguments play a fundamental role in the application of 
pseudoholomorphic curves to problems in symplectic and contact geometry: in the
closed case we have Gromov's compactness theorem, and more generally
the compactness theorems of Symplectic
Field Theory \cite{SFTcompactness} 
for punctured holomorphic curves in noncompact symplectic cobordisms.
As a rule, the singularities of the compactified moduli space have positive
virtual codimension, which translates into algebraic invariants if 
transversality is achieved.  In general however, even if the moduli space
of smooth curves is regular, multiple covers can appear in the 
compactification and make transversality impossible without abstract
perturbations, thus presenting a large technical complication.  

The motivating idea of this paper is that by restricting to a certain
geometrically natural class of holomorphic curves in low dimensional
settings, one can use topological constraints to prevent the aforementioned 
analytical difficulties from arising---in fact the compactified moduli space 
turns out to have a miraculously nice structure.  Examples of this phenomenon 
have been seen previously in the compactness arguments of 
\cite{HWZ:foliations} and~\cite{Wendl:OTfol}, both of
which deal with stable finite energy foliations on contact $3$--manifolds.
Roughly speaking, a finite energy foliation on a contact manifold
$(M,\xi)$ is an
$\RR$--invariant collection of pseudoholomorphic curves in $\RR\times M$
which project to a foliation of $M$ outside some set of closed Reeb orbits.
The foliation is called \emph{stable} if it deforms smoothly under sufficiently
small perturbations of the data on~$M$; in particular, this requires that
every leaf be parametrized by an embedded holomorphic curve of index~$1$
or~$2$.  As is shown in \cite{Wendl:BP1}, the class of holomorphic curves we
consider here consists (in the positive index case) of precisely those curves
which can be used to form finite energy foliations.

\begin{figure}
\includegraphics{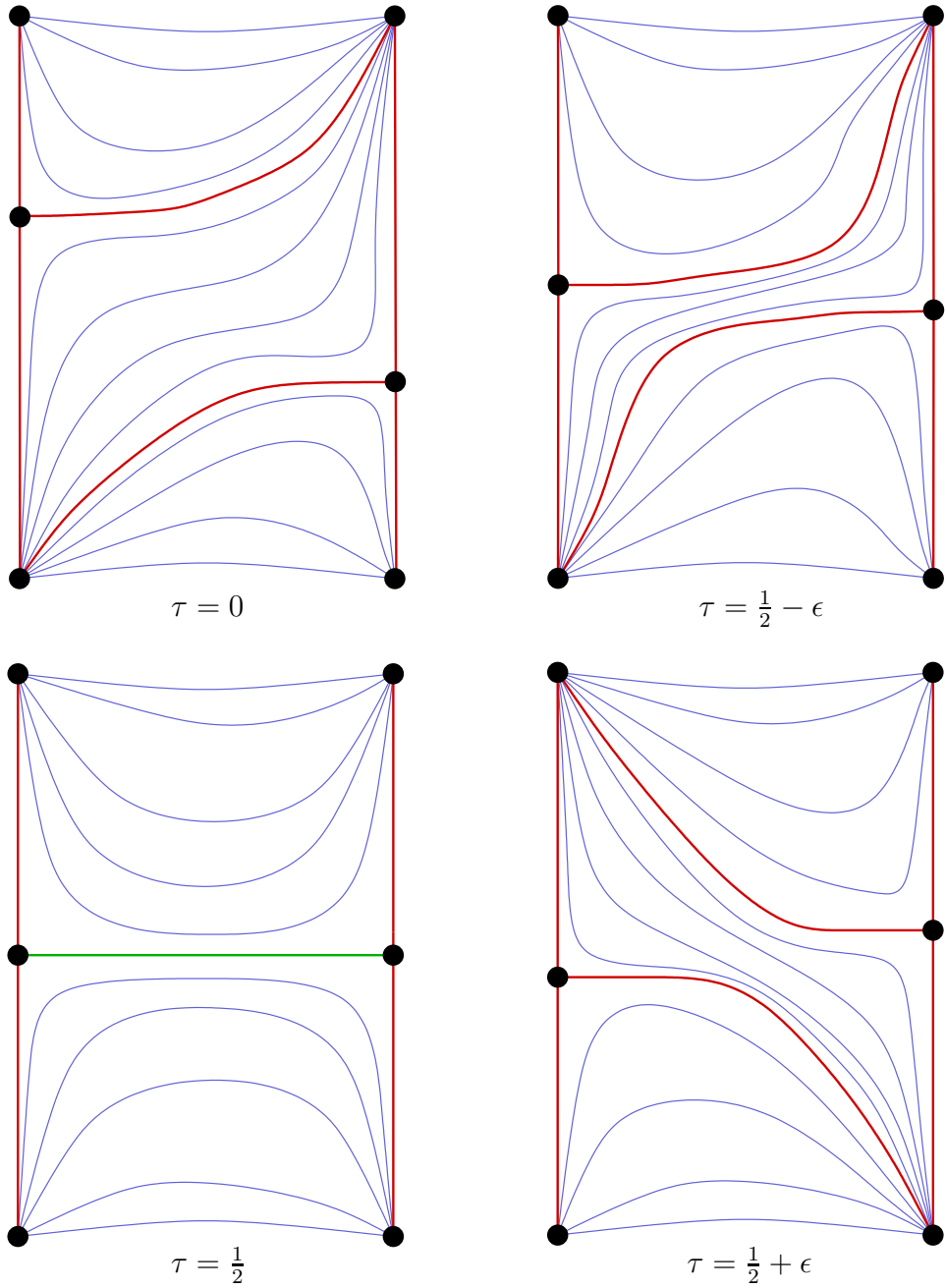}
\caption{\label{fig:fol_homotopy} 
Four steps in the deformation of a stable finite energy foliation 
under a generic homotopy of contact forms~$\{\lambda_\tau\}_{\tau\in[0,1]}$.
Each picture is a cross section, consisting mainly of
Reeb orbits that point through the page, index~$1$ holomorphic curves
(which appear isolated) and index~$2$ holomorphic curves (which appear in
$1$--parameter families).  The $\tau=1/2$ picture also contains a non-generic
index~$0$ holomorphic curve.}
\end{figure}

To illustrate the need for a compactness theorem, consider for a moment the
following question: can a stable finite energy foliation be deformed smoothly
under \emph{generic homotopies} of the contact form or complex structure?
Figure~\ref{fig:fol_homotopy} shows that the answer in general is no.  Here
we see a homotopy of the contact form which moves two
of the Reeb orbits that bound leaves of the foliation, and the families
of leaves
deform smoothly up until the isolated parameter value $\tau=1/2$.
At this value a non-generic index~$0$ leaf appears, producing
a discontinuous change in the structure of the
foliation.  The remarkable fact is that, at least in this example, 
the foliation \emph{survives} this discontinuous change: the leaves of the
unstable foliation at parameter $\tau=1/2$ can be glued to produce
a stable foliation for $\tau = 1/2 + \epsilon$.  To prove that this is what
should happen in general, one needs two fundamental ingredients:
\begin{itemize}
\item \emph{Compactness}: to show that the set of parameter values for
which a foliation exists is closed
\item \emph{Fredholm/gluing theory}: to show that that set is also open
\end{itemize}
The second ingredient 
only works if the linearized Cauchy-Riemann operator achieves
transversality: taking the homotopy to be sufficiently generic guarantees
this, but only for \emph{somewhere injective} holomorphic curves.  In this
regard, the standard compactness theory falls short, as it may in general
allow all manner of multiply covered curves to appear.  The result of this
paper is to strengthen the standard compactness theory accordingly for
the relevant class of holomorphic curves; this is a necessary step toward
carrying out the homotopy argument described above.

Along similar lines, M.~Hutchings \cite{Hutchings:index} has proved a
strong version of SFT-type compactness for a class of embedded index~$1$
and~$2$ curves in $4$--dimensional symplectizations, a result which forms
the analytical basis of Periodic Floer Homology and Embedded Contact
Homology.  The result proved here is different in several respects.
The condition on our set of curves is seemingly stricter than that of
Hutchings (though technically, neither implies the other), and the result is 
correspondingly
stronger: where Hutchings' limits allow certain types of multiple covers
(over trivial cylinders), ours do not.  In a different sense, the setup for
our main result is more general because it uses no genericity 
assumptions and is valid for arbitrary (also negative) Fredholm 
indices.  The restriction on multiple covers in the limit arises from 
topological considerations, independent of analysis; in particular
we make crucial use of the recently developed intersection theory for
punctured holomorphic curves, due to R.~Siefring \cite{Siefring:intersection}.

Hutchings' results suggest another possible application for our
compactness theory: it may be
possible to define specifically low-dimensional symplectic or contact
invariants (as in Gromov-Witten or Symplectic Field Theory \cite{SFT}) 
by counting 
this restricted class of holomorphic curves.  If such a theory exists, 
it has an immediate technical
advantage over general SFT, in that it seemingly can be defined without any
need for restrictive topological assumptions (e.g.~semipositivity)
or abstract perturbations.

The present work is part of a larger program involving compactness for a
special class of embedded holomorphic curves in $4$--dimensional symplectic
cobordisms.  We focus here on the special case where the 
target space is the $\RR$--invariant symplectization of
a $3$--manifold~$M$.  The relevant class of holomorphic curves is then 
distinguished by the property of being not only embedded in $\RR\times M$
but also having embedded projections to~$M$.  We'll give the required
definitions and state simple versions of the main theorems in 
\S\ref{subsec:mainresult}; these are implied by some slightly more
technical results which we state and prove in \S\ref{sec:mainproof},
after developing the necessary machinery.  
We will also give some more details in \S\ref{subsec:discussion}
on the general program into which this work fits, and state some partial
results for nontrivial symplectic cobordisms.

\subsection*{Acknowledgements}
This work benefited greatly from lengthy discussions with Michael
Hutchings and Richard Siefring, and I'd also like to thank Denis Auroux
and Cliff Taubes for some helpful suggestions.  Part of the work was
conducted during a visit to the FIM at ETH Z\"urich, and I thank Dietmar 
Salamon and the FIM for their hospitality.  Thanks also to the referee
for attentive reading and several suggestions that have improved the 
exposition.

\subsection{Setup and main results}
\label{subsec:mainresult}

The following structure was introduced in \cite{SFTcompactness} as a
general setting in which one has compactness results for punctured
holomorphic curves.
Let $M$ be a closed, oriented $3$--manifold.  We define a 
\emph{stable Hamiltonian structure} on $M$ to be a tuple
$\hH = (\xi,X,\omega,J)$, where\footnote{The tuple
($\xi,X,\omega)$, not including $J$, is equivalent to a \emph{framed stable 
Hamiltonian structure} in the definition given by \cite{Siefring:asymptotics}.
A similar definition appears in \cite{EliashbergKimPolterovich} with the
additional requirement that $\omega$ be exact.  The inclusion of $J$ in the 
data is not so natural geometrically, but convenient for our purposes.}
\begin{itemize}
\item $\xi$ is a smooth cooriented $2$--plane distribution on $M$
\item $\omega$ is a smooth closed $2$--form on $M$ which restricts to a 
 symplectic structure on the vector bundle $\xi \to M$
\item $X$ is a smooth vector field which is transverse to $\xi$, satisfies
 $\omega(X,\cdot) \equiv 0$, and whose flow preserves $\xi$
\item $J$ is a smooth complex structure on the bundle $\xi \to M$, compatible
 with $\omega$ in the sense that $\omega(\cdot, J\cdot)$ defines a
 bundle metric
\end{itemize}
It follows from these definitions that the flow of $X$ also preserves the
symplectic structure defined by $\omega$ on $\xi$, and the special
$1$--form $\lambda$ associated to $\xi$ and $X$ by the conditions
$$
\lambda(X) \equiv 1,
\qquad
\ker\lambda \equiv \xi,
$$
satisfies $d\lambda(X,\cdot) \equiv 0$.

An important example of a stable Hamiltonian structure arises when 
$\lambda$ is a \emph{contact form} on~$M$: then $d\lambda$ defines a
symplectic structure on the contact structure $\xi := \ker\lambda$, so
if $X_\lambda$ is the corresponding
Reeb vector field and $J$ is any complex structure on $\xi$ compatible with
$d\lambda$, we obtain a stable Hamiltonian structure in the form
$(\xi,X_\lambda,d\lambda,J)$.  A few non-contact examples may be
found in \cite{SFTcompactness}, some of which have also
appeared in applications, e.g.~in \cite{EliashbergKimPolterovich} 
and~\cite{Wendl:OTfol}.

We shall denote \emph{periodic orbits} of $X$ by $\gamma = (x,T)$, where
$T > 0$ and $x : \RR\to M$ satisfies $\dot{x} = X(x)$ and $x(T) = x(0)$.
If $x, x' : \RR \to M$ differ only by $x(t) = x'(t + c)$ for some $c \in \RR$,
we regard these as the same orbit $\gamma = (x,T) = (x',T)$.
We say that $\gamma$ has \emph{covering number} $k \in \NN$ if $T = k\tau$, 
where $\tau > 0$ is the \emph{minimal period}, i.e.~the smallest number
$\tau > 0$ such that $x(\tau) = x(0)$.  An orbit with covering number~$1$
is called \emph{simply covered}.  The \emph{$k$--fold cover} of
$\gamma = (x,T)$ will be denoted by
$$
\gamma^k = (x,kT).
$$
We shall occasionally abuse notation and regard $\gamma$ as a subset of $M$;
it should always be remembered that the orbit itself is specified by both this
subset and the period.

The open $4$--manifold $\RR\times M$ is called the \emph{symplectization}
of $M$, and it has a natural $\RR$--invariant almost complex structure
$\tilde{J}$ associated to any stable Hamiltonian structure 
$\hH = (\xi,X,\omega,J)$.
This is defined by $\tilde{J}\p_a = X$ and $\tilde{J}v = Jv$
for $v \in \xi$, where $a$ denotes the coordinate on the $\RR$--factor
and $\p_a$ is the unit vector in the $\RR$--direction.  We then consider
pseudoholomorphic (or \emph{$\tilde{J}$--holomorphic}) curves
$$
\tilde{u} = (a,u) : (\dot{\Sigma},j) \to (\RR\times M, \tilde{J}),
$$
where $\dot{\Sigma} = \Sigma \setminus\Gamma$, $(\Sigma,j)$ is a closed 
Riemann surface, $\Gamma\subset \Sigma$ is a finite set of punctures,
and by definition $\tilde{u}$ satisfies the nonlinear Cauchy-Riemann equation
$T\tilde{u} \circ j = \tilde{J} \circ T\tilde{u}$.  It is convenient
to think of $(\dot{\Sigma},j)$ as a Riemann surface with cylindrical
ends, and we will sometimes refer to neighborhoods of the punctures as 
\emph{ends} of~$\dot{\Sigma}$.

The \emph{energy} of a punctured pseudoholomorphic curve 
$\tilde{u} = (a,u) : (\dot{\Sigma},j) \to
(\RR\times M,\tilde{J})$ is defined by
$$
E(\tilde{u}) = E_\omega(\tilde{u}) + E_\lambda(\tilde{u}),
$$
where
\begin{equation}
\label{eqn:omegaEnergy}
E_\omega(\tilde{u}) = \int_{\dot{\Sigma}} u^*\omega
\end{equation}
is the so-called \emph{$\omega$--energy}, and
$$
E_\lambda(\tilde{u}) = \sup_{\varphi\in\tT} \int_{\Sigma} 
\tilde{u}^*(d\varphi \wedge \lambda),
$$
with $\tT := \{ \varphi\in C^\infty(\RR,[0,1])\ |\ \varphi' \ge 0 \}$.
An easy computation shows that both integrands are nonnegative whenever
$\tilde{u}$ is $\tilde{J}$--holomorphic, and such a curve is constant if
and only if $E(\tilde{u}) = 0$.  When $\tilde{u}$ is proper, connected,
$\tilde{J}$--holomorphic and satisfies $E(\tilde{u}) < \infty$, 
we call it a \emph{finite energy surface}.
As shown in \cites{Hofer:weinstein,HWZ:props1}, finite energy
surfaces have \emph{asymptotically cylindrical} behavior at the punctures:
this means the map $\tilde{u} : \dot{\Sigma} \to \RR\times M$
approaches $\{\pm\infty\} \times \gamma_z$ at each puncture 
$z \in \Gamma$, where $\gamma_z$ is a (perhaps multiply covered) 
periodic orbit of $X$.
(See Prop.~\ref{prop:asymptotics} for a precise statement.)  The sign 
in this expression partitions $\Gamma$ into positive and negative punctures
$\Gamma = \Gamma^+ \cup \Gamma^-$.

\begin{defn}
\label{defn:trivialCylinder}
The \emph{trivial cylinder} over a periodic orbit $\gamma = (x,T)$ of $X$
is the finite energy surface with one positive and one negative puncture
given by
$$
\tilde{u} : \RR\times S^1 \to \RR\times M : (s,t) \mapsto (Ts,x(Tt)).
$$
\end{defn}

Let $\varphi_X^t$ denote the time--$t$ flow of $X$, and
recall that a periodic orbit $\gamma = (x,T)$ of $X$ is called 
\emph{nondegenerate} if the linearized time--$T$ return map
$d\varphi_X^T(x(0))|_{\xi_{x(0)}}$ does not have~$1$ in its spectrum.
Choosing a unitary trivialization $\Phi$ of $\xi$ along $x$, one can
associate to any nondegenerate orbit $\gamma$ its \emph{Conley-Zehnder index}
$\muCZ^\Phi(\gamma) \in \ZZ$.  The odd/even parity of $\muCZ^\Phi(\gamma)$ is 
independent of the choice $\Phi$, and we call the orbit \emph{odd} or
\emph{even} accordingly.  Dynamically, even orbits are always hyperbolic,
elliptic orbits are always odd, and there can also exist odd
hyperbolic orbits, whose double covers are then even.  The following 
piece of terminology is borrowed from Symplectic Field Theory \cite{SFT},
where the orbits in question are precisely those which must be excluded in
order to define coherent orientations.

\begin{defn}
\label{defn:badOrbit}
A \emph{bad} orbit of $X$ is an even periodic orbit which 
is a double cover of an odd hyperbolic orbit.
\end{defn}

A stable Hamiltonian structure $\hH = (\xi,X,\omega,J)$ will be called
\emph{nondegenerate} if all periodic orbits of $X$ are nondegenerate, and
we will say that a sequence $\hH_k = (\xi_k,X_k,\omega_k,J_k)$ \emph{converges}
to $\hH = (\xi,X,\omega,J)$ if each piece of the data converges in the
$C^\infty$--topology on $M$.  We shall be concerned mainly with the following
special class of holomorphic curves.

\begin{defn}
\label{defn:nicely}
A finite energy surface $\tilde{u} = (a,u) : \dot{\Sigma} \to \RR\times M$
will be called \emph{nicely embedded} if the map
$u : \dot{\Sigma} \to M$ is embedded.
\end{defn}

\begin{figure}
\includegraphics{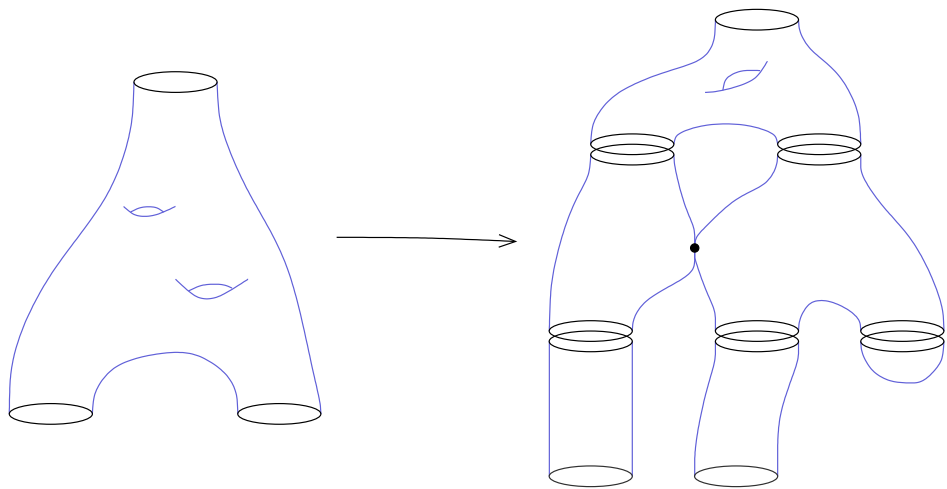}
\caption{\label{fig:building} A sequence of finite energy surfaces of
genus~$2$ converging to a stable holomorphic building
with three levels and arithmetic genus~$2$.  The middle level has a node.}
\end{figure}

By a compactness result in \cite{SFTcompactness}, sequences of finite 
energy surfaces with uniformly bounded genus and energy have subsequences
convergent to stable holomorphic \emph{buildings} (see 
Figure~\ref{fig:building}).  We will give precise
definitions in \S\ref{sec:buildings}; for now, let us simply recall
that a holomorphic building $\tilde{u}$ consists of finitely many
\emph{levels}, each of which is a (possibly disconnected) nodal
$\tilde{J}$--holomorphic curve with finite energy, and neighboring levels
can be attached to each other along matching \emph{breaking orbits}.
Every holomorphic building $\tilde{u}$ defines a graph $G_{\tilde{u}}$
whose vertices correspond to the smooth connected components of $\tilde{u}$,
with edges representing each node and breaking orbit.  We say that the
building $\tilde{u}$ is connected if the graph $G_{\tilde{u}}$ is
connected.

\begin{defn}
\label{defn:breakingOrbit}
For a holomorphic building $\tilde{u}$, a breaking orbit will be
called \emph{trivial} if deletion of the corresponding edge from
$G_{\tilde{u}}$ divides the graph into two components, one of which only
has vertices corresponding to trivial cylinders.  Breaking orbits
that do not have this property will be called \emph{nontrivial}.
\end{defn}

\begin{defn}
\label{defn:niceBuilding}
We say that a holomorphic building is \emph{nicely embedded} if
\begin{enumerate}
\item It has no nodes.
\item Each connected component is either a trivial cylinder or is 
nicely embedded.
\item If $\tilde{v}_1 = (b_1,v_1)$ and $\tilde{v}_2 = (b_2,v_2)$ are any two
distinct connected components, then the maps $v_1$ and $v_2$ 
either are identical or have no intersections.
\item Every nontrivial breaking orbit is even, and either simply
covered or both doubly covered and bad.
\end{enumerate}
\end{defn}

\begin{thm}
\label{thm:mainresult}
Assume $\hH_k = (\xi_k,X_k,\omega_k,J_k)$ is a sequence of stable Hamiltonian
structures converging to a nondegenerate stable Hamiltonian structure
$\hH = (\xi,X,\omega,J)$, 
and $\tilde{u}_k = (a_k,u_k) : \dot{\Sigma} \to \RR\times M$ are
$\tilde{J}_k$--holomorphic finite energy surfaces which are nicely embedded
and converge
in the sense of \cite{SFTcompactness} to a stable
$\tilde{J}$--holomorphic building
$\tilde{u}$.  Then $\tilde{u}$ is nicely embedded.
\end{thm}

Note that this statement assumes nothing about the index
of the curves $\tilde{u}_k$.  We will therefore obtain a stronger statement 
by restricting attention to generic data and curves of positive index.
Suppose $\tilde{u} : \dot{\Sigma} \to \RR\times M$ is a finite energy surface
with nondegenerate asymptotic orbits $\gamma_z$ at the 
punctures $z \in \Gamma$,
and $\Phi$ denotes a choice of unitary
trivialization for $\xi$ along each $\gamma_z$.
Then the Conley-Zehnder index of $\tilde{u}$ with respect to $\Phi$ is
defined to be the sum
\begin{equation}
\label{eqn:CZtotal}
\mu^\Phi(\tilde{u}) = \sum_{z\in\Gamma^+} \muCZ^\Phi(\gamma_z) -
\sum_{z\in\Gamma^-} \muCZ^\Phi(\gamma_z).
\end{equation}
Note that the parities of the orbits $\gamma_z$ partition $\Gamma$ into sets of
\emph{even} and \emph{odd} punctures, which we denote by
$$
\Gamma = \Gamma_0 \cup \Gamma_1.
$$
The \emph{Fredholm index} of $\tilde{u}$ is
\begin{equation}
\label{eqn:index}
\ind(\tilde{u}) = -\chi(\dot{\Sigma}) + 2c_1^\Phi(u^*\xi) + 
\muCZ^\Phi(\tilde{u}),
\end{equation}
where $c_1^\Phi(u^*\xi)$ denotes the relative first Chern number of the
bundle $u^*\xi \to \dot{\Sigma}$ with respect to $\Phi$, defined by
counting zeros of a generic section that is constant with respect to
$\Phi$ near the ends.
As shown in \cite{HWZ:props3}, $\ind(\tilde{u})$ is indeed the index of the 
linearized Cauchy-Riemann operator, and gives the virtual dimension of the 
moduli space of finite energy surfaces in a neighborhood of $\tilde{u}$.

For the following definition, it is useful to observe from \eqref{eqn:index}
that $\ind(\tilde{u}) + \Gamma_0$ is always even.
\begin{defn}
\label{defn:normalChern}
The \emph{normal first Chern number} $c_N(\tilde{u}) \in \ZZ$ of a finite 
energy surface $\tilde{u}$ of genus $g$ is defined by the relation
$$
2 c_N(\tilde{u}) = \ind(\tilde{u}) - 2 + 2g + \#\Gamma_0.
$$
\end{defn}
The meaning of $c_N(\tilde{u})$ is most easily seen by considering an
immersed, closed curve $\tilde{u} = (a,u) : \Sigma \to \RR\times M$: then
$2 c_N(\tilde{u}) = -\chi(\Sigma) + 2 c_1(u^*\xi) - 2 + 2g
= 2 c_1(u^*T(\RR\times M)) - 2\chi(\Sigma) = 2 c_1(N_{\tilde{u}})$, where
$N_{\tilde{u}} \to \Sigma$ is the normal bundle.  More generally, for immersed
curves with punctures, $c_N(\tilde{u})$ should be interpreted as
the relative first Chern number of $N_{\tilde{u}} \to \dot{\Sigma}$
with respect to special trivializations at the asymptotic orbits; this notion
will be made precise in \S\ref{sec:cN}.  The ``nicely embedded'' condition is
relevant to the normal first Chern number
for the following reason: $u$ is injective if and only if $\tilde{u}$ is
embedded and there is never
any intersection between $\tilde{u}$ and its $\RR$--translations
$\tilde{u}^c := (a + c,u)$ for $c \in \RR$.  In this case, $\tilde{u}$
belongs to a $1$--parameter family of non-intersecting finite energy surfaces,
which can be described via zero free sections of $N_{\tilde{u}}$.
This implies morally that $c_N(\tilde{u}) = 0$, a statement which becomes
literally true after applying the appropriate asymptotic constraints
(see~\S\ref{sec:constraints}).  In fact, one can show that for
generic $J$ (or generic parametrized families $J_\tau$), 
the condition $c_N(\tilde{u}) = 0$ is necessary
so that $\tilde{u}$ and all other finite energy surfaces nearby have
embedded projections into~$M$.  
A more detailed discussion of this may be found in \cite{Wendl:BP1}.
Note also that a linearized version of
positivity of intersections (see the discussion of $\windpi(\tilde{u})$
in \S\ref{sec:intersection}) implies $c_N(\tilde{u}) \ge 0$ for any
nicely embedded curve---thus $c_N(\tilde{u}) = 0$ is a minimality
condition.

According to Definition~\ref{defn:normalChern}, the condition 
$c_N(\tilde{u}) = 0$ allows exactly two cases where $\tilde{u}$ can have
positive index.  We will say that a nicely embedded curve
$\tilde{u}$ is \emph{stable} if either
\begin{itemize}
\item $\tilde{u}$ has index~$2$, genus~$0$ and no even punctures, or
\item $\tilde{u}$ has index~$1$, genus~$0$ and exactly one even puncture.
\end{itemize}

\begin{thm}
\label{thm:stable}
In addition to the assumptions of Theorem~\ref{thm:mainresult},
suppose the choice of $J$ in $\hH$ is generic and the curves 
$\tilde{u}_k$ are stable.  Then
\begin{enumerate}
\item
If $\ind(\tilde{u}) = 1$, $\tilde{u}$ is a stable nicely embedded finite 
energy surface, hence the moduli space of such curves up to
$\RR$--translation is compact.
\item
If $\ind(\tilde{u}) = 2$, then either $\tilde{u}$ is a nicely embedded
finite energy surface, or it is a building with exactly two nicely
embedded connected components, both stable with index~$1$, with projections
that do not intersect each other in $M$, and connected to
each other along a unique nontrivial breaking orbit.
\end{enumerate}
\end{thm}

\begin{figure}
\includegraphics{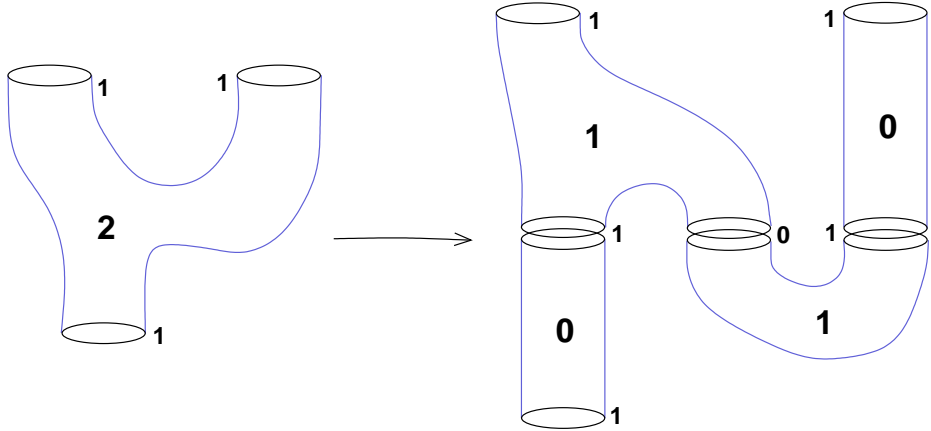}
\caption{\label{fig:stable} A sequence of stable nicely embedded finite
energy surfaces degenerating in accordance with Theorem~\ref{thm:stable}.
The numbers indicate the Fredholm indices of the components and parities 
of the orbits.
Note that each of the index~$0$ curves in the limit is a trivial cylinder,
and the odd breaking orbits are thus trivial. The limit has exactly
two nontrivial components (of index~$1$) and one nontrivial breaking 
orbit (even); the latter is also the unique even
orbit for each of the index~$1$ components.}
\end{figure}

\begin{figure}
\includegraphics{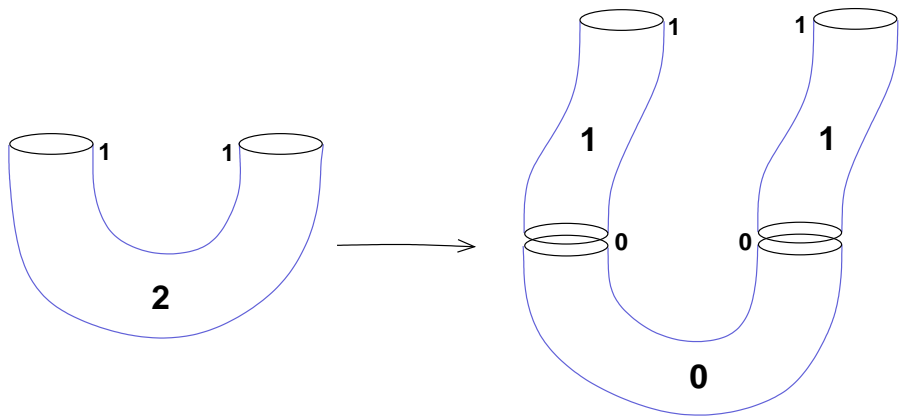}
\caption{\label{fig:nonGeneric} Convergence of stable nicely embedded
index~$2$ curves in the absence of genericity: the limit can now have more
than two nontrivial components, because nontrivial index~$0$ curves may
appear.  If the two odd orbits on the left are identical, it can also
occur that the two components of the top level in the limit are identical
curves, an outcome that is forbidden in the generic case.}
\end{figure}

Figure~\ref{fig:stable} 
shows a possible limit of stable index~$2$ curves.
Stranger things can happen if the genericity assumption is weakened:
for example if $J$ is not generic but belongs to a generic $1$--parameter
family $\{ J_\tau \}_{\tau \in \RR}$, then $\tilde{u}$ can contain
index~$0$ components (arbitrarily many, in principle) with
$\#\Gamma_0 = 2$, and there may be distinct nicely embedded components
with identical images (Figure~\ref{fig:nonGeneric}).

\subsection{Discussion}
\label{subsec:discussion}

The class of \emph{stable} nicely embedded finite energy surfaces that
we've defined above is of great interest in the theory of \emph{stable
finite energy foliations} introduced by Hofer, Wysocki and Zehnder
in \cite{HWZ:foliations}.  As shown in \cite{Wendl:BP1}, these are
precisely the curves whose moduli spaces form local foliations in
both $\RR\times M$ and~$M$.  Thus Theorem~\ref{thm:mainresult}, when combined
with some intersection theory and standard gluing analysis,
can be seen as a tool for proving stability of holomorphic foliations
under $\RR$--invariant homotopies.  

Likewise, Theorem~\ref{thm:stable}
guarantees a particularly nice structure for the moduli space of
leaves in a fixed foliation.  This provides the first half of the proof of 
an informal conjecture suggested in \cite{Wendl:thesis}, that to every 
stable finite energy foliation $\fF$ one can associate various 
SFT-type algebraic structures, in particular a Contact Homology
algebra $\HC_*(\fF)$.  In fact the result suggests more than
this, since it does not assume the existence of any foliation: one might
hope to encode this compactification of the space of nicely embedded 
index~$2$ curves
algebraically as in SFT, thus defining new invariants that
count nicely embedded index~$1$ curves.  In this case the
transversality problem is already solved for generic~$J$, so one would not
need any abstract perturbations or
restrictive assumptions on the target space.

To carry these ideas further, one needs
a corresponding compactness theorem for punctured holomorphic curves
$u : \dot{\Sigma} \to W$ in
nontrivial symplectic cobordisms $(W,\omega)$ with compatible~$J$.
In this case the ``nicely embedded'' condition makes no sense, but one
can formulate an appropriate generalization using the intersection
theory of punctured holomorphic curves 
defined in \cite{Siefring:intersection} (a conceptual 
summary without proofs may also be found in \cite{Wendl:automatic}).
In this theory, the standard adjunction formula for closed holomorphic
curves has a generalization of the form
$$
i(u,u) = 2\delta(u) + c_N(u) + \cov(u).
$$
Here $i(u,u)$ and $\delta(u)$ are generalizations of the homological
self intersection number and singularity number respectively: they are
homotopy invariant integers that count
intersections and singularities in addition to some nonnegative
``asymptotic terms'' (which vanish under generic perturbations).
The normal first Chern number $c_N(u)$ is again the integer given by
Definition~\ref{defn:normalChern}, and $\cov(u)$ is a nonnegative
integer that depends only on the asymptotic orbits of $u$:
it is zero if and only if all the relevant \emph{extremal eigenfunctions}
are simply covered (cf.~\S\ref{sec:eigenfunctions}).  Generically, one can now
characterize moduli spaces of nicely embedded curves in symplectizations by the
condition $i(\tilde{u},\tilde{u}) = 0$, and this is a sensible condition
to apply to certain curves in symplectic cobordisms as well.  Of
particular interest is the space of somewhere injective index~$2$ curves 
$u : \dot{\Sigma} \to W$ with $i(u,u) = 0$: these are automatically 
embedded and satisfy $c_N(u) = \cov(u) = 0$.  By a result in 
\cite{Wendl:BP1} in fact, such a curve is always regular and comes in a
smooth $2$--dimensional family of nonintersecting curves, which foliate
a neighborhood of $u(\dot{\Sigma})$ in $W$.

We can now state two partial results that we conjecture to be special 
cases of a more general theorem.
Assume for both that $(W,J)$ is an asymptotically cylindrical almost
complex manifold as in \cite{SFTcompactness}.

\begin{thmu}
Suppose $W$ is closed and $u_k : \Sigma \to W$ is a sequence of closed,
somewhere injective $J$--holomorphic curves with $\ind(u_k) = 2$ and
$i(u_k,u_k) = 0$, converging to a nodal curve~$u$.  Then $u$ is either
a smooth embedded curve or a nodal curve consisting of two embedded
index~$0$ components that intersect each other once transversely.  These
fit together with all smooth curves close to~$u$ as a singular foliation
of some neighborhood of the image of~$u$, with the nodal point as an
isolated singularity.
\end{thmu}
\begin{remarku}
We will not prove this here, but hope to include it in a future paper as a
special case of a much harder theorem for symplectic cobordisms.  The closed
case is comparatively simple and requires no substantially
new technology, only the
adjunction formula and some covering relations for $i(u,u)$ and $c_N(u)$.
It can also easily be generalized to apply to any 
$2$--dimensional moduli space of
curves that are embedded outside a set of marked points $z_1,\ldots,z_N$ 
satisfying fixed point constraints
$u(z_j) = p_j \in W$; one must then assume that $i(u,u)$ has the
smallest value allowed by the constraints.  The local structure of such
moduli spaces is studied in \cite{Wendl:BP1}, showing that locally they
form singular foliations in~$W$.  In this way one can also
accomodate immersed curves if the images of the self intersections are fixed.
\end{remarku}
\begin{thmu}[\cite{Wendl:automatic}]
Suppose $J$ is generic and $u_k : \dot{\Sigma} \to W$ is a sequence of 
embedded, punctured finite
energy $J$--holomorphic curves with $\ind(u_k) = 2$ and
$i(u_k,u_k) = 0$, and they converge to a smooth multiple cover
$u = v \circ \varphi$.  Then $v$ is an embedded index~$0$ curve with
$i(v,v) = -1$ and $u$ is immersed.  Moreover, the moduli space of curves
close to $u$ is a smooth orbifold, all other curves close to~$u$ are embedded, 
and they fit together with
$v$ as a foliation on some neighborhood of the image of~$v$.
\end{thmu}

In both cases, as with Theorems~\ref{thm:mainresult} and~\ref{thm:stable},
the upshot is that the degeneration in the limit is
nice enough so that transversality can still be achieved---this is
true even in the second case, despite the appearance of a multiple cover
in the limit.  (The latter can happen only in symplectic cobordisms that
are both noncompact and nontrivial).  
The reason one obtains smoothness in this case has to do with
the transversality results of Hofer, Lizan and Sikorav 
\cite{HoferLizanSikorav}, which are generalized in \cite{Wendl:automatic}:
specifically in dimension~$4$, one can sometimes use topological constraints
to prove transversality for all~$J$ (not just generic choices).  This
does not depend on $u$ being somewhere injective, though it is important
that $u$ is \emph{immersed}, and in fact the proof of the latter fact is also
based partly on such transversality arguments; 
see \cite{Wendl:automatic} for details.

\section{Holomorphic buildings in symplectizations}
\label{sec:buildings}

In this and the next few sections, we assemble some definitions and 
known results on 
punctured holomorphic curves and holomorphic buildings, fixing terminology 
and notation that will be used throughout.

Let $\dD$ denote the open unit disk in $\CC$, and write
$\dot{\dD} = \dD \setminus \{0\}$.  We define the 
\emph{circle compactification} of $\dot{\dD}$ as follows.
Using the biholomorphic map
$$
\varphi : (0,\infty)\times S^1 \to \dD\setminus\{0\} :
(s,t) \mapsto e^{- 2\pi(s+it)}
$$
to identify $\dot{\dD}$ with the half-cylinder, define 
$\overline{\dD} := \dot{\dD} \cup (\{\infty\}\times S^1) \cong 
(0,\infty]\times S^1$.  This is a 
topological surface with boundary, and has natural smooth structures over 
the interior $\interior{\overline{\dD}} = \dot{\dD}$ as well as the boundary 
$\p\overline{\dD} = \delta_0 := \{ \infty\}\times S^1$.

We use this to define a circle compactification $\overline{\Sigma}$ for 
$\dot{\Sigma} = \Sigma\setminus\Gamma$, where $(\Sigma,j)$ is any Riemann 
surface with isolated punctures $\Gamma\subset\Sigma$.
For each $z \in \Gamma$,
choose coordinates to identify a neighborhood of $z$ 
biholomorphically with $\dD$, identify the punctured neighborhood as above
with a half-cylinder and then add a \emph{circle at infinity}
$\delta_z \cong \{\infty\}\times S^1$ by replacing the half-cylinder with
$(0,\infty]\times S^1$.  The result is an oriented topological 
surface with boundary,
$$
\overline{\Sigma} = \dot{\Sigma} \cup 
\left( \bigcup_{z\in\Gamma} \delta_z\right),
$$
where the subsets $\interior{\overline{\Sigma}} = \dot{\Sigma}$ and
$\p\overline{\Sigma} = \bigcup_{z\in\Gamma} \delta_z$ inherit natural smooth 
structures that are independent of the choices of holomorphic coordinates.
The interior also has a conformal structure, and the complex structure on
$T_z\Sigma$ for $z \in \Gamma$ defines a special class of 
diffeomorphisms $\varphi : S^1 \to \delta_z$, which are all related to each 
other by a constant shift, i.e.~$\varphi_1(t) = \varphi_2(t + \text{const})$.  
For any two punctures $z_1,z_2 \in \Gamma$, an orientation reversing
diffeomorphism $\psi : \delta_{z_1} \to \delta_{z_2}$ will be called 
\emph{orthogonal} if it can be written as $\psi(t) = -t$ with respect to
some choice of special diffeomorphisms $\delta_{z_i} \cong S^1$.
Observe that  $\overline{\Sigma}$ is compact if $\Sigma$ is closed.

A closed \emph{nodal Riemann surface} with marked points consists of the data
$$
\mathbf{S} = (S,j,\Gamma,\Delta),
$$
where $(S,j)$ is a closed (but not necessarily connected) Riemann surface,
and $\Gamma, \Delta \subset S$ are disjoint finite subsets with the 
following additional structure:
\begin{itemize}
\item $\Gamma$ is ordered,
\item elements of $\Delta$ are grouped into pairs
$\overline{z}_1,\underline{z}_1,\ldots,\overline{z}_n,\underline{z}_n$.
\end{itemize}
We call $\Delta$ the \emph{double points} of $\mathbf{S}$, and
$\Gamma$ the \emph{marked points}.
Let $\dot{S} = S \setminus (\Gamma \cup \Delta)$, with circle
compactification $\overline{S}$.  For a pair
$\{\overline{z},\underline{z}\} \subset \Delta$, a \emph{decoration} at
$\{\overline{z},\underline{z}\}$ is an orientation reversing orthogonal
diffeomorphism $\psi : \delta_{\overline{z}} \to \delta_{\underline{z}}$, and a
decoration $\psi$ of $\mathbf{S}$ is a choice of decorations at all
pairs $\{\overline{z},\underline{z}\} \subset \Delta$; we can regard this as a
diffeomorphism on a certain subset of $\p\overline{S}$.
We call $\mathbf{S} := (S,j,\Gamma,\Delta,\psi)$ a
\emph{decorated nodal Riemann surface}.

Given $\mathbf{S}$ with decoration $\psi$, define
$$
\overline{\mathbf{S}} = \overline{S} 
/ \{ z \sim \psi(z) \}.
$$
This is an oriented topological surface with boundary, with a conformal
structure that degenerates at $\p\overline{\mathbf{S}} =
\bigcup_{z\in\Gamma} \delta_z$
and also at a certain set of disjoint circles $\Theta_\Delta \subset
\overline{\mathbf{S}}$, one for each double point pair
$\{\overline{z},\underline{z}\} \subset \Delta$.
There is a natural inclusion of $\dot{S}$ into
$\overline{\mathbf{S}}$ as the subset
$$
\dot{S} = \interior{\overline{\mathbf{S}}} \setminus
\Theta_\Delta.
$$

We say that the nodal surface $\mathbf{S}$ is \emph{connected} if
$\overline{\mathbf{S}}$ is connected, and define its \emph{arithmetic genus} 
to be the genus of $\overline{\mathbf{S}}$.  Neither of these definitions
depends on the choice of decoration.

Let $M$ be a closed $3$--manifold with stable Hamiltonian structure
$\hH = (\xi,X,\omega,J)$ and associated almost complex structure $\tilde{J}$.
If $\mathbf{S} = (S,j,\Gamma,\Delta)$ is a closed nodal Riemann surface 
with marked points, a \emph{nodal $\tilde{J}$--holomorphic curve}
$$
\tilde{u} : \mathbf{S} \to \RR\times M
$$
is a proper finite energy 
pseudoholomorphic map $\tilde{u} = (a,u) : (S \setminus\Gamma,j) \to
(\RR\times M,\tilde{J})$ 
such that for each pair $\{\overline{z},\underline{z}\}
\subset \Delta$, $\tilde{u}(\overline{z}) = \tilde{u}(\underline{z})$.
In this context each pair $\{ \overline{z},\underline{z} \} \subset \Delta$ is
called a \emph{nodal pair}, or simply a \emph{node} of
$\tilde{u}$.  The marked points $\Gamma$ are called \emph{punctures} of
$\tilde{u}$, and the asymptotic behavior of $a : S\setminus\Gamma \to \RR$
determines the \emph{sign} of each, defining a partition
$\Gamma = \Gamma^+ \cup \Gamma^-$.
Observe that for any decoration $\psi$ of $\mathbf{S}$,
$\tilde{u} : S\setminus\Gamma \to M$ has a natural continuous extension
$$
(\bar{a},\bar{u}) : \overline{\mathbf{S}} \to [-\infty,\infty] \times M,
$$
which is constant on each connected component of $\Theta_\Delta$ and
maps $\p\overline{\mathbf{S}}$ to $\{\pm\infty\} \times M$;
in particular the restriction of $\bar{u}$ 
to each $\delta_z \subset \p\overline{\mathbf{S}}$ for
$z \in \Gamma^\pm$ defines a positively/negatively oriented parametrization 
of a periodic orbit $\gamma_z$ of $X$.

Consider next a collection of nodal $\tilde{J}$--holomorphic curves
$$
\tilde{u}_m = (a_m,u_m) : \mathbf{S}_m = (S_m, j_m,
\Gamma_m,\Delta_m) \to \RR\times M
$$
for $m=1,\ldots,n$.
Denote $\p_\pm \overline{\mathbf{S}}_m = \bigcup_{z\in\Gamma_m^\pm}
\delta_z$, and suppose there are orientation reversing orthogonal 
diffeomorphisms
$$
\varphi_m : \p_+\overline{\mathbf{S}}_m \to
\p_-\overline{\mathbf{S}}_{m+1}
$$
for each $m=1,\ldots,n-1$.  Then the collection
$$
\tilde{u} = (\tilde{u}_1,\ldots,\tilde{u}_n ;
\varphi_1,\ldots,\varphi_{n-1})
$$
is called a \emph{$\tilde{J}$--holomorphic building of height~$n$} if
for each $m=1,\ldots,n-1$,
$$
\bar{u}_m|_{\p_+\overline{\mathbf{S}}_m} = 
\bar{u}_{m+1} \circ \varphi_m.
$$
The nodal curves $\tilde{u}_m$ are called \emph{levels} of
$\tilde{u}$.  For each $m=1,\ldots,n-1$ and $\underline{z} \in \Gamma_m^+$,
there is a unique $\overline{z} \in \Gamma_{m+1}^-$ such that
$\varphi_m(\delta_{\underline{z}}) = \delta_{\overline{z}}$.  We then call the
pair $\{\overline{z},\underline{z}\}$ a \emph{breaking pair}, and denote
by $\gamma_{(\overline{z},\underline{z})}$ the \emph{breaking orbit}
parametrized by
$\bar{u}_m|_{\delta_{\underline{z}}}$ and 
$\bar{u}_{m+1}|_{\delta_{\overline{z}}}$.  
Let $\Delta_{\oC}$ denote the set of all punctures in
$\Gamma_1 \cup \ldots \cup \Gamma_n$ that belong to breaking pairs.

Define the \emph{partially decorated}
nodal Riemann surface $\mathbf{S} = (S,j,\Gamma,\Delta,\varphi)$, 
where $(S,j)$ is the disjoint union of $(S_1,j_1),\ldots,
(S_n,j_n)$, $\Gamma = \Gamma^+ \cup \Gamma^- :=
\Gamma_n^+ \cup \Gamma_1^-$, $\Delta$ is the union of the
breaking pairs in $\Delta_{\oC}$ with the nodal pairs in
$\Delta_{\oN} := \Delta_1 \cup \ldots \cup \Delta_n$, and
$\varphi$ is the collection of decorations at
the breaking pairs $\{ \overline{z},\underline{z} \}$ defined by
$\varphi_m : \p_+\overline{\mathbf{S}}_m \to
\p_-\overline{\mathbf{S}}_{m+1}$.  We will call $\mathbf{S}$ the
\emph{domain} of $\tilde{u}$, and indicate this via the shorthand notation
$$
\tilde{u} : \mathbf{S} \to \RR\times M.
$$
Choosing arbitrary decorations
$\psi_m$ for each $\mathbf{S}_m$, these together with $\varphi$
define a decoration $\psi$ for $\mathbf{S}$,
and $\overline{\mathbf{S}}$ is now the surface obtained from
$\overline{\mathbf{S}}_1 \cup \ldots \cup
\overline{\mathbf{S}}_n$ by gluing boundaries together via $\varphi$.
There is then a continuous map
$$
\bar{u} : \overline{\mathbf{S}} \to M
$$
such that $\bar{u}|_{\overline{\mathbf{S}}_m} = \bar{u}_m$.
The orbits $\gamma_z$ parametrized by $\bar{u}|_{\delta_z}$ for 
$z \in \Gamma^\pm$ are called \emph{asymptotic orbits} of~$\tilde{u}$.

\begin{remark}
Technically, what we've defined should be called holomorphic buildings
\emph{with zero marked points}, since all the marked points of $\mathbf{S}$
are being viewed as punctures of $\tilde{u}$.  One can also define
holomorphic buildings with marked points, though we will not
need them here; see \cite{SFTcompactness} for details.
\end{remark}

The relationship of a building $\tilde{u}$ with its domain $\mathbf{S}$
gives rise to a slightly more general notion which we will find useful.
\begin{defn}
Suppose $\mathbf{S} = (S,j,\Gamma,\Delta,\psi)$ is a nodal Riemann surface
with $\Delta$ partitioned into two sets $\Delta_{\operatorname{C}} 
\cup \Delta_{\operatorname{N}}$, each organized in pairs,
called \emph{breaking} pairs and \emph{nodal} pairs respectively, 
and $\psi$ denotes a choice of decoration at each of the breaking pairs.
A \emph{generalized $\tilde{J}$--holomorphic building} $\tilde{u} :
\mathbf{S} \to \RR\times M$ is then a proper
finite energy $\tilde{J}$--holomorphic
map $\tilde{u} = (a,u) : S\setminus (\Gamma \cup \Delta_{\oC}) \to \RR\times M$ 
such that
\begin{enumerate}
\item For each nodal pair $\{ \overline{z},\underline{z} \} \subset \Delta_{\oN}$,
$\tilde{u}(\overline{z}) = \tilde{u}(\underline{z})$.
\item Completing $\psi$ to a decoration of $\mathbf{S}$ by choosing 
arbitrary decorations at the nodal pairs, 
$u : S\setminus (\Gamma \cup \Delta)$ extends to a continuous map
$\bar{u} : \overline{\mathbf{S}} \to M$.
\end{enumerate}
\end{defn}
Considering orientations, we see that each breaking pair
$\{ \overline{z},\underline{z} \} \subset \Delta_{\oC}$ includes one positive
and one negative puncture of $\tilde{u}$.

Just as with nodal Riemann surfaces, we say that the generalized 
building $\tilde{u} :
\mathbf{S} \to \RR\times M$ is \emph{connected} if $\overline{\mathbf{S}}$ is
connected, and its \emph{arithmetic genus} is defined as the genus of
$\overline{\mathbf{S}}$.  A \emph{connected
component} of $\tilde{u}$ is the finite energy surface obtained by restricting
the map $\tilde{u} : S\setminus(\Gamma \cup \Delta_{\oC}) 
\to \RR\times M$ to any connected component of its domain.
The sets $\Gamma^\pm \subset S$ are the positive
and negative \emph{punctures} of $\tilde{u}$.  In general, each connected
component may have some punctures that do not belong to
$\Gamma$ but are included among the breaking pairs $\Delta_{\oC}$: we call
these \emph{breaking punctures}.

Every holomorphic building is also a generalized holomorphic building in an
obvious way.  The main difference is that the components of generalized
buildings cannot in general be assigned \emph{levels}, and every component
may have positive and negative punctures which are not breaking punctures;
holomorphic buildings have these only at the top and bottom levels.

\begin{defn}
A generalized building $\tilde{u} : \mathbf{S} \to \RR\times M$ is
\emph{stable} if each connected component $\dot{S}_i \subset S \setminus
(\Gamma \cup \Delta)$ on which $\tilde{u}$ is constant
satisfies $\chi(\dot{S}_i) < 0$.
\end{defn}

\begin{notation}
For any generalized holomorphic building $\tilde{u}$ and puncture 
$z \in \Gamma$, we will always denote by
$$
\gamma_z := \bar{u}(\delta_z)
$$
the asymptotic orbit at $z$, and for breaking pairs $\{ \overline{z},
\underline{z} \} \subset \Delta_{\oC}$ denote the breaking orbit
$\bar{u}(\delta_{\overline{z}}) = \bar{u}(\delta_{\underline{z}})$ by
$$
\gamma_{(\overline{z},\underline{z})} = \gamma_{\overline{z}} =
\gamma_{\underline{z}}.
$$
Unless stated otherwise, the domain $\mathbf{S}$ will be assumed to
consist of the data $(S,j,\Gamma,\Delta,\varphi)$ as defined above.
When multiple domains are under discussion, we'll often use $\mathbf{S}'$
to denote a second domain $(S',j',\Gamma',\Delta',\varphi')$, or similarly
$\mathbf{S}\trivial = (S\trivial,j\trivial,\Gamma\trivial,\Delta\trivial,
\varphi\trivial)$ and so forth.
\end{notation}

\begin{defn}
\label{defn:augmentation}
Given a generalized holomorphic building
$\tilde{u} : \mathbf{S} \to \RR\times M$, an \emph{augmentation}
of $\tilde{u}$ at a puncture $z \in \Gamma^\pm$ is the generalized building
$\tilde{u}' : \mathbf{S}' \to \RR\times M$, defined as follows:
\begin{enumerate}
\item
$(S',j')$ is the disjoint union of $(S,j)$ with a sphere $(S_T,j_T) :=
(S^2,i)$.  Denote by $p_-,p_+ \in S'$ the points $0$ and $\infty$ respectively
in $S_T$.
\item
$\Gamma' = (\Gamma \cup \{p_\pm\} ) \setminus \{z\}$.
\item
$\Delta'_{\oC}$ is $\Delta_{\oC}$ with the addition of one extra pair
$\{z,p_\mp\}$.
\item
$\Delta'_{\oN} = \Delta_{\oN}$.
\item
$\tilde{u}'|_{S\setminus(\Gamma \cup \Delta_{\oC})} = \tilde{u}$ and
$\tilde{u}'|_{S_T \setminus \{p_+,p_-\}}$ is a trivial cylinder over
$\gamma_z$.  A decoration is then chosen at $\Delta'_{\oC} 
\setminus \Delta_{\oC}$ so 
that $\tilde{u}'$ is a generalized holomorphic building.
\end{enumerate}
An augmentation at a breaking pair $\{\overline{z},\underline{z}\}
\subset \Delta_{\oC}$ is defined in the same manner with the following changes:
\begin{enumerate}
\item
$\Gamma' = \Gamma$.
\item
$\Delta'_{\oC}$ is $\Delta_{\oC}$ with two additional pairs,
$\{ \underline{z},p_- \}$ and $\{ p_+,\overline{z} \}$.
\item
$\tilde{u}'|_{S_T \setminus \{p_+,p_-\}}$ is a trivial cylinder over
$\gamma_{(\overline{z},\underline{z})}$.
\end{enumerate}
In general, an \emph{augmentation} of $\tilde{u}$ is any generalized building
obtained from $\tilde{u}$ by a finite sequence of these two operations.
\end{defn}
Augmentation is essentially the operation of shifting levels of $\tilde{u}$
by inserting trivial cylinders.  One should think of $\tilde{u}'$ as being
\emph{homotopic} to $\tilde{u}$, generalizing the fact that a finite energy 
surface $\tilde{v} = (b,v)$ is homotopic to any of its $\RR$--translations
$\tilde{v}^c = (b + c,v)$ for $c \in \RR$; an augmentation is in some sense
a sequence of infinite $\RR$--translations.

Just as one can insert trivial cylinders in a generalized building, one can
also ``collapse'' them.
\begin{defn}
\label{defn:core}
Suppose $\tilde{u} : \mathbf{S} \to \RR\times M$ is a generalized 
holomorphic building such that at least
one connected component is not a trivial cylinder.
The \emph{core} of $\tilde{u}$ is then the unique
generalized building $\tilde{u}\core : \mathbf{S}\core \to
\RR\times M$ such that $\tilde{u}$ is an augmentation
of $\tilde{u}\core$ and no connected component of 
$\tilde{u}\core$ is a trivial cylinder.
\end{defn}

\begin{defn}
\label{defn:subbuilding}
Given a generalized building $\tilde{u} : \mathbf{S} \to \RR\times M$,
a \emph{subbuilding} $\tilde{u}' : \mathbf{S}' \to \RR\times M$ of
$\tilde{u}$ is a generalized building such that
\begin{enumerate}
\item 
$S'$ is an open and closed subset of $S$, on which $j' = j$.
\item
$\Gamma'$ is the union of $\Gamma \cap S'$ with all
$z \in S'$ for which $\{z,z'\}$ is a breaking pair in $\Delta_{\oC}$ with
$z' \not\in S'$.
\item
$\Delta'_{\oC}$ is the set of all 
breaking pairs $\{\overline{z},\underline{z}\}$ 
in $\Delta_{\oC}$ for which both $\overline{z}$ and $\underline{z}$ are in $S'$,
and $\psi'$ is the restriction of $\psi$.
\item
$\Delta'_{\oN}$ is the set of all nodal pairs $\{\overline{z},\underline{z}\}$
in $\Delta_{\oN}$ for which both $\overline{z}$ and $\underline{z}$ are in $S'$.
\item
$\tilde{u}' = \tilde{u}|_{S' \setminus (\Gamma'\cup \Delta'_{\oC})}$
\end{enumerate}
\end{defn}

We refer to \cite{SFTcompactness} for the detailed definition of
what it means for a sequence of finite energy surfaces to converge to a 
holomorphic building.  We will only need to use the following fact,
immediate from the definition:
\begin{prop}
\label{prop:convergence}
If $\tilde{u}_k = (a_k,u_k) : (\dot{\Sigma}_k,j_k) \to 
(\RR\times M,\tilde{J}_k)$ is a sequence of 
finite energy surfaces converging to a $\tilde{J}$--holomorphic building 
$\tilde{u} : \mathbf{S} \to \RR\times M$, then for sufficiently large~$k$ 
there exist homeomorphisms
$\varphi_k : \overline{\mathbf{S}} \to \overline{\Sigma}_k$, 
restricting to smooth maps $\overline{\mathbf{S}}\setminus
(\p\overline{\mathbf{S}} \cup \Theta_\Delta) \to \dot{\Sigma}_k$,
such that
\begin{enumerate}
\item
$u_k \circ \varphi_k \to u$ in $\Cinftyloc(\overline{\mathbf{S}} \setminus
(\p\overline{\mathbf{S}} \cup \Theta_\Delta),M)$,
\item
$\bar{u}_k \circ \varphi_k \to \bar{u}$ 
in $C^0(\overline{\mathbf{S}},M)$.
\end{enumerate}
\end{prop}

The proof of Theorem~\ref{thm:mainresult} will rely heavily on our ability
to control the normal first Chern number for components of a holomorphic
building.  Positivity of intersections guarantees that this number is
nonnegative for any finite energy surface that is not preserved by the
$\RR$--action.  The hardest part of the proof therefore involves the
so-called \emph{trivial curves}, for which this number can be negative.

\begin{defn}
\label{defn:trivialCurve}
A finite energy surface $\tilde{u} : \dot{\Sigma} \to \RR\times M$ will be
called a \emph{trivial curve} if $E_\omega(\tilde{u}) = 0$, and
\emph{nontrivial} if $E_\omega(\tilde{u}) > 0$.
\end{defn}

Examining the integrand in \eqref{eqn:omegaEnergy}, one finds that a
finite energy surface $\tilde{u} = (a,u) : \dot{\Sigma} \to \RR\times M$ is
trivial if and only if the image of $du(z)$ is everywhere tangent to $X$,
which means $u(\dot{\Sigma})$ is contained in a single periodic orbit
$\gamma$.  If $\tilde{u}$ is not constant, then this implies 
one can always write $\tilde{u} = \tilde{v} \circ
\varphi$ where
$\varphi : \dot{\Sigma} \to \RR\times S^1$ is a holomorphic branched cover
and $\tilde{v}$ is the trivial cylinder over $\gamma$.

\begin{defn}
\label{defn:trivialBuilding}
A holomorphic building or generalized holomorphic building 
will be called \emph{trivial} if it is connected, has no nodes, and all its 
connected components are trivial curves.  Additionally, such a building
will be called \emph{cylindrical} if it has arithmetic genus zero and
exactly two punctures.
\end{defn}

Observe that every nonconstant trivial curve has at least one positive and one
negative puncture, and the same statement therefore holds for trivial
(generalized) buildings.  It follows that in general,
a trivial building $\tilde{u} : \mathbf{S} \to \RR\times M$ has
$\chi(\overline{\mathbf{S}}) \le 0$, with equality if and only if $\tilde{u}$
is cylindrical.

\begin{prop}
\label{prop:trivialCurves}
A nonconstant trivial building is cylindrical if and only if it is an
augmentation of a trivial cylinder.
\end{prop}
\begin{proof}
Observe first that the statement is true for any building $\tilde{u}$
with only one level (i.e.~a finite energy surface),
for then $\tilde{u}$ covers a trivial cylinder $\tilde{v}$ by a holomorphic
map $\varphi : \RR\times S^1 \to \RR\times S^1$ of degree $k \in \NN$, and
every such map is of the form $\varphi(s,t) = (ks,kt)$ up to constant
shifts in $s$ and~$t$.

For a more general trivial building $\tilde{u} : \mathbf{S} \to \RR\times M$
with two punctures and $\chi(\overline{\mathbf{S}}) = 0$, we only need
observe that under these assumptions, no connected component of $\tilde{u}$
can have nontrivial genus or more than two punctures.
\end{proof}

\section{Asymptotic eigenfunctions}
\label{sec:eigenfunctions}

Let $\gamma = (x,T)$ be a periodic orbit of $X$, and writing $S^1 :=
\RR / \ZZ$, define the parametrization
$$
\mathbf{x} : S^1 \to M : t \mapsto x(Tt).
$$
We can then view the normal bundle to $\gamma$ as the induced bundle
$\mathbf{x}^*\xi \to S^1$.  Choosing any symmetric connection $\nabla$
on $M$, we define the \emph{asymptotic operator}
$$
\mathbf{A}_\gamma : \Gamma(\mathbf{x}^*\xi) \to \Gamma(\mathbf{x}^*\xi) :
v \mapsto - J (\nabla_t v - T \nabla_v X).
$$
One can check that this expression doesn't depend on the choice
$\nabla$ and gives a well defined section of $\mathbf{x}^*\xi$.  Morally,
$\mathbf{A}_\gamma$ is the Hessian of a certain action functional 
on $C^\infty(S^1,M)$, whose critical points are the closed characteristics
of~$X$.  As an unbounded operator on $L^2(\mathbf{x}^*\xi)$ with domain
$H^1(\mathbf{x}^*\xi)$, $\mathbf{A}_\gamma$ is self adjoint, with spectrum
$\sigma(\mathbf{A}_\gamma)$ 
consisting of isolated real eigenvalues of multiplicity 
at most two.  We shall sometimes refer to the eigenfunctions of
$\mathbf{A}_\gamma$ as \emph{asymptotic eigenfunctions}.
Recall that the orbit $\gamma$ is degenerate if and only if $0 \in 
\sigma(\mathbf{A}_\gamma)$.

Operators of this form are fundamental in the asymptotic
analysis of punctured holomorphic curves, as demonstrated by the following 
result proved in \cites{HWZ:props1,Mora,Siefring:asymptotics}.  Denote
$$
\RR_+ = [0,\infty), \quad
\RR_- = (-\infty,0], \quad
Z_\pm = \RR_\pm \times S^1,
$$
and assign to $Z_\pm$ the standard complex structure
$i \frac{\p}{\p s} = \frac{\p}{\p t}$ in terms of the
coordinates $(s,t) \in Z_\pm$.
We will use the term \emph{asymptotically constant reparametrization} to mean
a smooth embedding $\varphi : Z_\pm \to Z_\pm$ for which there are constants
$s_0 \in \RR$ and $t_0 \in S^1$ such that $\varphi(s - s_0,t - t_0) - (s,t) 
\to 0$ with all derivatives as $s \to \pm\infty$.

\begin{prop}
\label{prop:asymptotics}
Suppose $\tilde{u} = (a,u) : Z_\pm \to \RR_\pm\times M$ is a proper
finite energy half-cylinder asymptotic to a 
nondegenerate orbit $\gamma = (x,T)$.  
Then there is an asymptotically constant reparametrization
$\varphi : Z_\pm \to Z_\pm$ such that for $|s|$ sufficiently large,
$\tilde{u}\circ\varphi :
Z_\pm \to \RR\times M$ is either $(Ts,x(Tt))$ or satisfies the following 
asymptotic formula:
there exists an eigenfunction $e \in \Gamma(\mathbf{x}^*\xi)$ of 
$\mathbf{A}_\gamma$ with negative/positive
eigenvalue $\lambda$ such that
$$
\tilde{u} \circ \varphi(s,t) = \exp_{(Ts,x(Tt))}\left[ e^{\lambda s} \cdot
\left( e(t) + r(s,t) \right) \right],
$$
where $\exp$ is defined with respect to any  $\RR$--invariant connection 
on~$\RR\times M$ and
$r(s,t) \in \xi_{x(Tt)}$ satisfies $r(s,t) \to 0$ with all derivatives as
$s \to \pm\infty$.
\end{prop}

Since nontrivial eigenfunctions $e \in \Gamma(\mathbf{x}^*\xi)$
are never zero,
this implies in particular that $u(s,t)$ is either contained in $\gamma$
or never intersects $\gamma$ for sufficiently large $|s|$.  Similar
formulas are proved in \cite{Siefring:asymptotics} for the relative asymptotic 
behavior of distinct holomorphic half-cylinders approaching the same orbit; 
this is fundamental to the intersection theoretic results that we will review
in \S\ref{sec:intersection}.

Prop.~\ref{prop:asymptotics} obviously applies to finite energy surfaces
$\tilde{u} : \dot{\Sigma} \to \RR\times M$ by identifying a punctured 
disk-like neighborhood of each puncture $z \in \Gamma^\pm$
biholomorphically with $Z_\pm$.

\begin{defn}
\label{defn:convergenceRate}
For $\tilde{u} : Z_\pm \to \RR_\pm\times M$ as in Prop.~\ref{prop:asymptotics},
if $\tilde{u}$ satisfies the asymptotic formula with a nontrivial eigenfunction
$e$ and eigenvalue $\lambda$, we call $|\lambda| > 0$ the \emph{transversal
convergence rate} of $\tilde{u}$ and say that $e$ \emph{controls} the
asymptotic approach of $\tilde{u}$ to $\gamma$.  Otherwise, if $\tilde{u}$
is simply a reparametrization of $(Ts,x(Tt))$ near infinity, we define the
transversal convergence rate to be $+\infty$.  Similar wording will be used
also for more general punctured holomorphic curves 
$\tilde{u} : \dot{\Sigma} \to
\RR\times M$ approaching orbits $\gamma_z$ at $z \in \Gamma$.
\end{defn}

It will be important to understand the eigenfunctions of $\mathbf{A}_\gamma$ 
in greater detail.  For $k \in \NN$, define a parametrization of $\gamma^k$ by
$$
\mathbf{x}_k : S^1 \to M : t \mapsto x(kTt).
$$
Choose a unitary trivialization $\Phi$ of 
$\mathbf{x}^*\xi$, and use $\Phi$ also to denote the natural trivialization
induced on $\mathbf{x}_k^*\xi$.
Every nowhere zero section $v \in \Gamma(\mathbf{x}^*\xi)$ now has a well 
defined \emph{winding number}
$$
\wind^\Phi(v) \in \ZZ.
$$
By a result in \cite{HWZ:props2}, the winding number of a nontrivial
eigenfunction $e$ of $\mathbf{A}_\gamma$ depends only on its 
eigenvalue $\lambda$, thus we can sensibly write $\wind^\Phi(e) = 
\wind^\Phi(\lambda)$.
In fact, the result in question proves:
\begin{prop}
\label{prop:HWZwinding}
$\wind^\Phi : \sigma(\mathbf{A}) \to \ZZ$ is a monotone increasing 
function, and
for each $k \in \ZZ$, there are precisely two eigenvalues $\lambda \in
\sigma(\mathbf{A}_\gamma)$ counted with multiplicity 
such that $\wind^\Phi(\lambda) = k$.
\end{prop}
Assuming $\gamma$ to be nondegenerate, we now define the integers
\begin{equation}
\label{eqn:defAlpha}
\begin{split}
\alpha^\Phi_-(\gamma) &= \max \{ \wind^\Phi(\lambda) \ |\ 
\text{$\lambda \in \sigma(\mathbf{A}_\gamma)$, $\lambda < 0$} \}, \\
\alpha^\Phi_+(\gamma)  &= \min \{ \wind^\Phi(\lambda) \ |\ 
\text{$\lambda \in \sigma(\mathbf{A}_\gamma)$, $\lambda > 0$} \}, \\
p(\gamma)      &= \alpha^\Phi_+(\gamma) - \alpha^\Phi_-(\gamma),
\end{split}
\end{equation}
noting that the \emph{parity} $p(\gamma) \in \{0,1\}$ doesn't depend on $\Phi$.  
Another result in
\cite{HWZ:props2} then gives the following formula for the Conley-Zehnder
index:
\begin{equation}
\label{eqn:CZwinding}
\muCZ^\Phi(\gamma) = 2\alpha^\Phi_-(\gamma) + p(\gamma) = 2\alpha^\Phi_+(\gamma) - p(\gamma).
\end{equation}

\begin{defn}
We will say that a nontrivial eigenfunction $e$ of $\mathbf{A}_\gamma$ is
a \emph{positive/negative extremal eigenfunction} of $\gamma$ if
$\wind^\Phi(e) = \alpha^\Phi_\pm(\gamma)$.
\end{defn}

If $e \in \Gamma(\mathbf{x}^*\xi)$ satisfies
$\mathbf{A}_\gamma e = \lambda e$, we define the $k$--fold cover
$e^k \in \Gamma(\mathbf{x}_k^*\xi)$
by $e^k(t) = e(kt)$ and find $\mathbf{A}_{\gamma^k} e^k = k\lambda e^k$.
In general, an eigenfunction $f$ of $\mathbf{A}_\gamma$ is called a
\emph{$k$--fold cover} if there is an orbit $\zeta$ and eigenfunction
$e$ of $\mathbf{A}_\zeta$ such that $\zeta^k = \gamma$ and $e^k = f$.
We say that $f$ is \emph{simply covered} if it is not a $k$--fold cover for
any $k > 1$.

\begin{lemma}
A nontrivial eigenfunction $f$ of $\mathbf{A}_{\gamma^k}$ is a $k$--fold
cover if and only if $\wind^\Phi(f) \in k \ZZ$.
\end{lemma}
\begin{proof}
Clearly if $f = e^k$ then $\wind^\Phi(f) = k\wind^\Phi(e) \in k\ZZ$.
To see the converse, note that by Prop.~\ref{prop:HWZwinding} there is
a two-dimensional space of eigenfunctions $e$ of $\mathbf{A}_\gamma$ having
any given integer value of $\wind^\Phi(e)$.  This gives rise to a 
two-dimensional space of $k$--fold covers $e^k$ with winding
$k\wind^\Phi(e)$.  Since this attains all winding numbers in $k\ZZ$, every
eigenfunction of $\mathbf{A}_{\gamma^k}$ that is \emph{not} a $k$--fold cover
has winding in $\ZZ \setminus k\ZZ$.
\end{proof}

\begin{prop}
\label{prop:relPrime}
Let $\gamma$ be a simply covered periodic orbit, $\Phi$ a unitary 
trivialization of $\xi$ along $\gamma$ and $k \in \NN$.  Then a nontrivial
eigenfunction $e$ of 
$\mathbf{A}_{\gamma^k}$ is simply covered if and only if $k$ and $\wind^\Phi(e)$
are relatively prime.
\end{prop}
\begin{proof}
From the lemma, we see that $e$ is an $n$--fold cover if and only if 
$n$ divides both $k$ and $\wind^\Phi(e)$.  So $e$ is simply covered if and 
only if this is not true for any $n \in \{2,\ldots,k\}$.
\end{proof}

\section{Intersection theory}
\label{sec:intersection}

If $\tilde{u} = (a,u) : \dot{\Sigma} \to \RR\times M$ is a finite energy
surface and $\pi : TM \to \xi$ is the fiberwise linear projection along~$X$,
then the composition $\pi \circ Tu : T\dot{\Sigma} \to \xi$ defines a
section 
$$
\piTu : \dot{\Sigma} \to \Hom_\CC(T\dot{\Sigma},u^*\xi).
$$
It is shown in \cite{HWZ:props2} that $\piTu$
satisfies the similarity principle, thus it is
either trivial or has only isolated positive zeros; the latter is the case
unless $E_\omega(\tilde{u}) = 0$.  Assuming $E_\omega(\tilde{u}) > 0$ and
$\tilde{u}$ also has nondegenerate asymptotic orbits, the asymptotic 
formula of Prop.~\ref{prop:asymptotics} implies that $\piTu$ has no zeros 
outside some compact subset, and its winding near infinity is controlled by 
eigenfunctions of asymptotic operators.  We can thus define the integer
$$
\windpi(\tilde{u}) \ge 0
$$
as the algebraic count of zeros of $\piTu$.  It follows from the nonlinear
Cauchy-Riemann equation that $\windpi(\tilde{u}) = 0$
if and only if $u : \dot{\Sigma} \to M$ is immersed and transverse to~$X$.

Recalling the formula for $c_N(\tilde{u})$ from 
Definition~\ref{defn:normalChern}, a result in \cite{HWZ:props2}
shows that all finite energy surfaces $\tilde{u}$ with 
$E_\omega(\tilde{u}) > 0$ satisfy $\windpi(\tilde{u}) \le c_N(\tilde{u})$. 
Actually one can state this in a slightly stronger form as an equality.  For
$z \in \Gamma^\pm$, let $e_z$ be an asymptotic 
eigenfunction that controls the
approach of $\tilde{u}$ to $\gamma_z$, choose a unitary trivialization $\Phi$
of $\xi$ along $\gamma_z$ and define the \emph{asymptotic defect} 
of $\tilde{u}$ at $z$ to be the nonnegative integer
$$
\asympdef^z(\tilde{u}) = 
\left| \alpha^\Phi_\mp(\gamma_z) - \wind^\Phi(e_z) \right|.
$$
This is zero if and only if the asymptotic approach to $\gamma_z$ is
controlled by an extremal eigenfunction.
The total asymptotic defect of $\tilde{u}$ is then defined as
$$
\asympdef(\tilde{u}) = \sum_{z\in\Gamma} \asympdef^z(\tilde{u}).
$$
Now the argument in \cite{HWZ:props2} implies:

\begin{prop}
\label{prop:windpi}
For any finite energy surface $\tilde{u}$ with $E_\omega(\tilde{u}) > 0$
and nondegenerate asymptotic orbits,
$$
\windpi(\tilde{u}) + \asympdef(\tilde{u}) = c_N(\tilde{u}),
$$
and both terms on the left hand side are nonnegative.
\end{prop}

Next we collect some important results from the intersection theory 
of finite energy surfaces, due to R.~Siefring.
In the following, we assume all periodic orbits of $X$ are nondegenerate
and write punctured holomorphic disks as half-cylinders
$Z_\pm \to \RR_\pm\times M$.  These statements, proved in 
\cite{Siefring:intersection}, are all based on the construction
of homotopy invariant ``asymptotic intersection numbers'', which are
well defined due to the relative asymptotic formulas proved 
in~\cite{Siefring:asymptotics}.

\begin{prop}
\label{prop:uEmbedded}
Suppose $\tilde{u} = (a,u) : Z_\pm \to \RR_\pm\times M$ is a proper
finite energy half-cylinder such that $u : Z_\pm \to M$ is embedded.
Then any asymptotic eigenfunction controlling $\tilde{u}$ at infinity is
simply covered.
\end{prop}

\begin{prop}
\label{prop:uvDisjoint}
Suppose $\tilde{u} = (a,u) : Z_\pm \to \RR_\pm\times M$ and
$\tilde{v} = (b,v) : Z_\pm \to \RR_\pm\times M$ are proper finite 
energy half-cylinders
asymptotic to $\gamma^m$ and $\gamma^n$ respectively for some simply 
covered orbit $\gamma$
and $m,n \in \NN$.  Assume also that $u$ and $v$ are both embedded and do not
intersect each other.  Then $m = n$, and the asymptotic 
eigenfunctions controlling $\tilde{u}$ and $\tilde{v}$ at infinity have the
same winding number.  
\end{prop}

\begin{prop}
\label{prop:oppositeSigns}
Suppose $\tilde{u}_+ : (a_+,u_+) : Z_+ \to \RR_+ \times M$ and
$\tilde{u}_- : (a_-,u_-) : Z_- \to \RR_- \times M$ are proper finite energy 
half-cylinders asymptotic to $\gamma^{k_+}$ and $\gamma^{k_-}$ respectively 
for some simply covered periodic orbit $\gamma$ and $k_\pm \in \NN$.  
Assume also $u_+$ and $u_-$ are both embedded and do not intersect each other.
Then both have asymptotic defect zero, and either
\begin{enumerate}
\setlength{\itemsep}{0in}
\item $\gamma$ is even and $k_+ = k_- = 1$, or
\item $\gamma$ is odd hyperbolic and $k_+ = k_- = 2$, 
hence $\gamma^{k_+} = \gamma^{k_-}$ is bad.
\end{enumerate}
\end{prop}

\section{Constraints at the asymptotic orbits}
\label{sec:constraints}

Given a periodic orbit $\gamma$, a positive/negative 
\emph{asymptotic constraint} for $\gamma$ is a real number
$c \ge 0$ such that $\mp c \not\in\sigma(\mathbf{A}_\gamma)$.
We will say that a proper finite energy half-cylinder $\tilde{u} : Z_\pm \to
\RR_\pm \times M$ asymptotic to $\gamma$ is \emph{compatible} with this
constraint if its transversal convergence rate (recall
Definition~\ref{defn:convergenceRate}) is strictly greater than~$c$.
Similarly, for a generalized holomorphic building $\tilde{u} :
\mathbf{S} \to \RR\times M$ with punctures $\Gamma = \Gamma^+ \cup \Gamma^-$,
denote by $\mathbf{c} = \{ c_z \}_{z\in\Gamma}$ an association of a 
positive/negative asymptotic
constraint $c_z$ to each asymptotic orbit $\gamma_z$ for $z \in \Gamma^\pm$, 
and say that $\tilde{u}$ is compatible with $\mathbf{c}$ if for every 
$z \in \Gamma$, the corresponding end has transversal convergence rate 
strictly greater than~$c_z$.  Observe that the space of holomorphic buildings 
compatible with a given set of asymptotic constraints is a closed subset 
of the space of all holomorphic buildings.

Let $\gamma$ be a nondegenerate orbit and fix a unitary trivialization
$\Phi$ of $\xi$ along $\gamma$.  Then if $c$ is a positive asymptotic
constraint for $\gamma$, define the positive
\emph{constrained Conley-Zehnder index} by
\begin{equation}
\label{eqn:CZconstrained1}
\muCZ^\Phi(\gamma ; c) = \muCZ^\Phi(\gamma) - 
\# \left( \sigma(\mathbf{A}_\gamma) \cap (-c,0)\right),
\end{equation}
where eigenvalues are counted with multiplicity.  Similarly, for
$c$ a negative asymptotic constraint, define the negative
constrained Conley-Zehnder index
\begin{equation}
\label{eqn:CZconstrained2}
\muCZ^\Phi(\gamma ; -c) = \muCZ^\Phi(\gamma) +
\# \left( \sigma(\mathbf{A}_\gamma) \cap (0,c)\right),
\end{equation}
and if $\tilde{u} : \mathbf{S} \to \RR\times M$ is a holomorphic
building compatible with constraints $\mathbf{c}$, choose unitary 
trivializations $\Phi$ for $\xi$ along all asymptotic orbits $\gamma_z$ and
define the total constrained Conley-Zehnder index
$$
\muCZ^\Phi(\tilde{u} ; \mathbf{c}) =
\sum_{z \in \Gamma^+} \muCZ^\Phi(\gamma_z ; c_z) -
\sum_{z \in \Gamma^-} \muCZ^\Phi(\gamma_z ; -c_z).
$$
The even/odd parity of the constrained indices $\muCZ^\Phi(\gamma_z ; \pm c_z)$
for $z \in \Gamma^\pm$ defines a \emph{constrained parity} for each 
puncture, thus defining a new partition
$$
\Gamma = \Gamma_0(\mathbf{c}) \cup \Gamma_1(\mathbf{c}).
$$
Now define the \emph{constrained Fredholm index}
\begin{equation}
\label{eqn:indexConstrained}
\ind(\tilde{u} ; \mathbf{c}) =
- \chi(\overline{\mathbf{S}}) + 2 c_1^\Phi(\bar{u}^*\xi) +
\muCZ(\tilde{u} ; \mathbf{c}).
\end{equation}
As shown in \cite{Wendl:BP1} (based on arguments in \cite{HWZ:props3}),
if $\tilde{u}$ is a finite energy surface,
$\ind(\tilde{u} ; \mathbf{c})$ is the virtual dimension of the
moduli space of finite energy surfaces near $\tilde{u}$ that are compatible
with~$\mathbf{c}$.

The relation between Conley-Zehnder indices and winding numbers has a
straightforward generalization to the constrained case.
Given an orbit $\gamma$ and $c \in \RR$ with $-c \not\in 
\sigma(\mathbf{A}_\gamma)$, define
\begin{equation}
\label{eqn:defAlphaConstrained}
\begin{split}
\alpha^\Phi_-(\gamma; c) &= \max \{ \wind^\Phi(\lambda) \ |\ 
\text{$\lambda \in \sigma(\mathbf{A}_\gamma)$, $\lambda < -c$} \}, \\
\alpha^\Phi_+(\gamma; c)  &= \min \{ \wind^\Phi(\lambda) \ |\ 
\text{$\lambda \in \sigma(\mathbf{A}_\gamma)$, $\lambda > -c$} \}, \\
p(\gamma;c)      &= \alpha^\Phi_+(\gamma;c) - 
\alpha^\Phi_-(\gamma;c).
\end{split}
\end{equation}
Then combining \eqref{eqn:CZwinding} with \eqref{eqn:CZconstrained1}
and~\eqref{eqn:CZconstrained2}, we have
\begin{equation}
\label{eqn:CZwindingConstrained}
\muCZ^\Phi(\gamma;c) = 2\alpha^\Phi_-(\gamma;c) + p(\gamma;c) 
= 2\alpha^\Phi_+(\gamma;c) - p(\gamma;c).
\end{equation}
A nontrivial eigenfunction $e$ of $\mathbf{A}_\gamma$ will now be called a
positive/negative \emph{extremal} eigenfunction \emph{with respect to the
constraint $|c|$} if $\wind^\Phi(e) = \alpha^\Phi_\pm(\gamma ; c)$.

Now if $\tilde{u} : \dot{\Sigma} \to \RR\times M$ is a finite energy surface
compatible with $\mathbf{c}$ and $E_\omega(\tilde{u}) > 0$, define the
\emph{constrained asymptotic defect} at $z \in \Gamma^\pm$ by
$$
\asympdef^z(\tilde{u} ; c_z) =
\left| \alpha^\Phi_\mp(\gamma_z ; \pm c_z) - \wind^\Phi(e_z) \right|,
$$
where $e_z$ is an eigenfunction controlling the asymptotic approach of
$\tilde{u}$ to~$\gamma_z$.  The total \emph{constrained asymptotic defect}
is then
$$
\asympdef(\tilde{u} ; \mathbf{c}) := \sum_{z \in \Gamma}
\asympdef^z(\tilde{u} ; c_z).
$$
This sum is nonnegative, and is zero if and only if $\tilde{u}$ is
controlled by extremal eigenfunctions with respect to the constraints
at every puncture.

If $\tilde{u}$ is a generalized building
compatible with constraints $\mathbf{c}$, then every subbuilding
$\tilde{u}_0$ is compatible with a natural
set of \emph{induced constraints} $\hat{\mathbf{c}}$ defined as follows.
For each puncture $z$ of $\tilde{u}_0$ that is also a puncture of
$\tilde{u}$, set $\hat{c}_z = c_z$, and for all other punctures of
$\tilde{u}_0$ (i.e.~those which are only \emph{breaking}
punctures of $\tilde{u}$), set $\hat{c}_z = 0$.  The following relation is
easily verified using \eqref{eqn:indexConstrained}.

\begin{prop}
\label{prop:indexAdditivity}
If $\tilde{u} : \mathbf{S} \to \RR\times M$ is a generalized holomorphic
building compatible with constraints $\mathbf{c}$ and it has connected
components $\tilde{u}_i$ with induced constraints $\mathbf{c}_i$, then
$$
\ind(\tilde{u} ; \mathbf{c}) = \sum_i \ind(\tilde{u}_i ; \mathbf{c}_i) +
\#\Delta_{\operatorname{N}}.
$$
\end{prop}

Suppose now that $\tilde{u}'$ is an augmentation of $\tilde{u}$: there is then
a canonical bijection between the sets of punctures for each, so
a set of asymptotic constraints $\mathbf{c}$ on either induces one on the
other, which we'll also denote by $\mathbf{c}$.  
However if $\tilde{u}'$ is compatible with $\mathbf{c}$, it is
\emph{not} necessarily true that $\tilde{u}$ is as well.  Indeed, it may
happen that for a given puncture $z \in \Gamma$, the component of 
$\tilde{u}'$ containing $z$ is a trivial cylinder, and is therefore
compatible with arbitrarily strict asymptotic constraints, which is not
necessarily true for $\tilde{u}$.  On the other hand, if $\tilde{u}'$
arises as the limit of a sequence $\tilde{u}_k$ of finite energy surfaces
compatible with $\mathbf{c}$, then in a neighborhood of $z \in \Gamma$,
convergence to $\tilde{u}'$ and convergence to $\tilde{u}$ are
equivalent notions.  It follows that both $\tilde{u}$ and $\tilde{u}'$
are in this case compatible with $\mathbf{c}$: this is true in particular
if $\tilde{u}$ is the core of $\tilde{u}'$.

\section{The normal first Chern number}
\label{sec:cN}

Suppose $\tilde{u} : \mathbf{S} \to \RR\times M$ is a generalized 
holomorphic building compatible with asymptotic constraints $\mathbf{c} = 
\{c_z\}_{z\in\Gamma}$, and the asymptotic orbits $\gamma_z$ are 
all nondegenerate.  Fix a unitary trivialization $\Phi$ for $\xi$ along
each~$\gamma_z$.

\begin{defn}
\label{defn:C0normalChern}
Define the \emph{constrained normal first Chern number} of $\tilde{u}$
with respect to $\mathbf{c}$ as the integer
$$
c_N(\tilde{u} ; \mathbf{c}) = c_1^\Phi(\bar{u}^*\xi) - 
\chi(\overline{\mathbf{S}}) + 
\sum_{z\in \Gamma^+} \alpha^\Phi_-(\gamma_z ; c_z) -
\sum_{z\in \Gamma^-} \alpha^\Phi_+(\gamma_z ; -c_z).
$$
\end{defn}
One can easily check that this doesn't depend on $\Phi$, and a simple
computation using \eqref{eqn:CZwindingConstrained} and 
\eqref{eqn:indexConstrained} shows that
\begin{equation}
\label{eqn:cNindex}
2 c_N(\tilde{u} ; \mathbf{c}) = 
\ind(\tilde{u} ; \mathbf{c}) - 2 + 2g + \#\Gamma_0(\mathbf{c}),
\end{equation}
where $g$ is the arithmetic genus of $\tilde{u}$.
The new formula is therefore consistent with Definition~\ref{defn:normalChern}.

The following result is immediate from the definition.
\begin{prop}
If $\tilde{u}'$ is an augmentation of $\tilde{u}$ then
$c_N(\tilde{u}' ; \mathbf{c}) = c_N(\tilde{u} ; \mathbf{c})$.
\end{prop}

We also have an immediate generalization of Prop.~\ref{prop:windpi}:
\begin{prop}
\label{prop:windpiConstrained}
If $\tilde{u}$ is a finite energy surface compatible with $\mathbf{c}$
and $E_\omega(\tilde{u}) > 0$, then
$$
\windpi(\tilde{u}) + \asympdef(\tilde{u} ; \mathbf{c}) = c_N(\tilde{u} ;
\mathbf{c}),
$$
and both terms on the left hand side are nonnegative.
\end{prop}

\begin{prop}
\label{prop:cNsum}
Suppose $\tilde{u} : \mathbf{S} \to \RR\times M$ is a generalized
holomorphic building compatible with asymptotic constraints $\mathbf{c}$,
and $\tilde{u}_i : \dot{S}_i
\to \RR\times M$ are the connected components of $\tilde{u}$, with
induced constraints $\mathbf{c}_i$ for $i=1,\ldots,N$.  Then
$$
c_N(\tilde{u} ; \mathbf{c}) = \sum_{i=1}^N c_N(\tilde{u}_i ; \mathbf{c}_i) + 
\sum_{\{\overline{z},\underline{z}\}\subset \Delta_{\oC}} 
p(\gamma_{(\overline{z},\underline{z})})
+ \#\Delta_{\oN}.
$$
\end{prop}
\begin{proof}
We must check that $c_N$ behaves appropriately under certain
natural operations on generalized holomorphic buildings.
The simplest such operation is the \emph{disjoint union} of two buildings
$\tilde{u} : \mathbf{S} \to \RR\times M$ and
$\tilde{u}' : \mathbf{S}' \to \RR\times M$ with constraints $\mathbf{c}$ and
$\mathbf{c}'$ respectively: this defines a building
$\tilde{u} \sqcup \tilde{u}' :
\mathbf{S} \sqcup \mathbf{S}' \to \RR\times M$, compatible with the
obvious union of constraints $\mathbf{c} \sqcup \mathbf{c}'$.  Clearly then,
$$
c_N(\tilde{u} \sqcup \tilde{u}' ; \mathbf{c} \sqcup \mathbf{c}') = 
c_N(\tilde{u} ; \mathbf{c}) + c_N(\tilde{u}' ; \mathbf{c}').
$$

Next, if $\tilde{u} : \mathbf{S} \to \RR\times M$ is a building with
two points $z,z' \in S\setminus (\Gamma \cup \Delta)$
such that $\tilde{u}(z) = \tilde{u}(z')$, we can add a node to $\tilde{u}$
and define
$\odot_{(z,z')} \tilde{u} : \odot_{(z,z')} \mathbf{S} \to \RR\times M$ 
by adding $\{z,z'\}$ to the set of nodal pairs.  
This decreases the Euler characteristic of $\overline{\mathbf{S}}$ by~$2$,
thus
$$
c_N(\odot_{(z,z')} \tilde{u} ; \mathbf{c}) =
c_N(\tilde{u} ; \mathbf{c}) + 2.
$$
Similarly, if there are punctures $\underline{z} \in \Gamma^+$ and
$\overline{z} \in \Gamma^-$ for which $\gamma_{\underline{z}} =
\gamma_{\overline{z}}$ and $c_{\underline{z}} = c_{\overline{z}} = 0$, 
then we can change $\tilde{u}$
by ``gluing'' these punctures, which means adding $\{\overline{z},
\underline{z}\}$ to the set of breaking pairs and choosing an
appropriate decoration so that the result is a generalized holomorphic building
$\boxplus_{(\overline{z},\underline{z})} \tilde{u} :
\boxplus_{(\overline{z},\underline{z})}\mathbf{S} \to \RR\times M$.
By losing two unconstrained punctures, this operation
subtracts $\alpha^\Phi_-(\gamma_{\underline{z}}) -
\alpha^\Phi_+(\gamma_{\overline{z}}) = 
-p(\gamma_{(\overline{z},\underline{z})})$ from $c_N(\tilde{u} ; \mathbf{c})$, 
hence
$$
c_N(\boxplus_{(\overline{z},\underline{z})} \tilde{u} ; \mathbf{c})
= c_N(\tilde{u} ; \mathbf{c}) + p(\gamma_{(\overline{z},\underline{z})}).
$$
Composing these operations as often as necessary gives the stated result.
\end{proof}

\section{Proofs of the main results}
\label{sec:mainproof}

We will now state and prove stronger, more technical versions of
Theorems~\ref{thm:mainresult} and~\ref{thm:stable}.
Assume $\hH_k = (\xi_k,X_k,\omega_k,J_k)$
is a sequence of stable Hamiltonian structures converging to
$\hH = (\xi,X,\omega,J)$, where the latter is nondegenerate,
and $\tilde{u} : \mathbf{S} \to \RR\times M$ is a $\tilde{J}$--holomorphic
building compatible 
with asymptotic constraints $\mathbf{c} = \{ c_z \}_{z\in\Gamma}$.
Since the orbits $\gamma_z$ are nondegenerate, for sufficiently large~$k$ 
there are unique periodic orbits $\gamma_{z,k}$ of $X_k$ such that
$$
\gamma_{z,k} \to \gamma_z,
$$
in the sense that these orbits have parametrizations $S^1 \to M$ that
converge in the $C^\infty$--topology.
We may also assume that for each $z \in \Gamma^\pm$,
$\mp c_z \not\in\sigma(\mathbf{A}_{\gamma_{z,k}})$.

\begin{thm}
\label{thm:mainTechnical}
Assume $\tilde{u}_k$ are nicely embedded $\tilde{J}_k$--holomorphic finite 
energy surfaces converging to the $\tilde{J}$--holomorphic 
building $\tilde{u}$,
such that the $\tilde{u}_k$ are also compatible with $\mathbf{c}$ and 
$c_N(\tilde{u}_k ; \mathbf{c}) = 0$.  Then
$\tilde{u}$ is nicely embedded, its core $\tilde{u}\core$ 
is compatible with $\mathbf{c}$, and for every connected component 
$\tilde{v}_i$ of $\tilde{u}\core$ with induced constraints $\mathbf{c}_i$,
$c_N(\tilde{v}_i ; \mathbf{c}_i) = 0$.
\end{thm}

\begin{remark}
If $\tilde{u}_k \to \tilde{u}$ under the assumptions of 
Theorem~\ref{thm:mainresult}, then one can assume after taking a subsequence
that all the $\tilde{u}_k$ are compatible with some choice of asymptotic
constraints $\mathbf{c}$ such that $c_N(\tilde{u}_k ; \mathbf{c}) =
\windpi(\tilde{u}_k) + \asympdef(\tilde{u}_k ; \mathbf{c}) = 0$.
This is why Theorem~\ref{thm:mainTechnical} implies 
Theorem~\ref{thm:mainresult}.  Similarly, Theorem~\ref{thm:stable} is
a special case of the next statement.
\end{remark}

\begin{thm}
\label{thm:stableTechnical}
In addition to the assumptions of Theorem~\ref{thm:mainTechnical},
suppose $J$ is generic.  Then $\ind(\tilde{u} ; \mathbf{c})$ is either~$1$
or~$2$.  If it is~$1$, then $\tilde{u}$ is a finite energy surface, 
hence the moduli space of such curves
with constraint~$\mathbf{c}$ is compact.  
If $\ind(\tilde{u} ; \mathbf{c}) = 2$ and $\tilde{u}$ is not a finite
energy surface, then it has exactly two nontrivial connected components
$\tilde{v}_i = (b_i,v_i)$, both with $\ind(\tilde{v}_i ; \mathbf{c}_i) = 1$,
such that $v_1$ and $v_2$ have no intersections in~$M$ and they are
connected to each other by a unique nontrivial breaking orbit.
\end{thm}

We begin now with some preparations for the proof of 
Theorem~\ref{thm:mainTechnical}.  By Prop.~\ref{prop:convergence},
we can assume without loss of generality that the curves $\tilde{u}_k$
have a fixed domain $\dot{\Sigma} = \Sigma\setminus\Gamma$ with varying
complex structures $j_k$, and there is a fixed homeomorphism
$$
\psi : \overline{\mathbf{S}} \to \overline{\Sigma}
$$
such that $u_k \circ \psi \to u$ in 
$\Cinftyloc(\overline{\mathbf{S}} \setminus (\p\overline{\mathbf{S}} \cup
\Theta_\Delta),M)$ and
$\bar{u}_k \circ \psi \to \bar{u}$ in $C^0(\overline{\mathbf{S}},M)$.
The punctures $\Gamma^\pm$ of $\tilde{u}_k$ and $\tilde{u}$ are also
identified via $\psi$, so we shall use the same notation for both: the
asymptotic orbit of $\tilde{u}_k$ at $z \in \Gamma$ is then
$\gamma_{z,k}$ for sufficiently large~$k$.

\begin{lemma}
\label{lemma:nontrivialCurve}
The building $\tilde{u}$ has at least one nontrivial component.
\end{lemma}
\begin{proof}
If $\tilde{u} : \mathbf{S} \to \RR\times M$ is a trivial
building, then $\bar{u} : \overline{\mathbf{S}} \to M$ represents the
trivial homology class $[\bar{u}] = 0 \in H_2(M)$.  Perturbing $\tilde{u}$
to $\tilde{u}_k$ with asymptotic orbits $\gamma_{z,k}$ for sufficiently
large~$k$, we also have $[\bar{u}_k] = 0 \in H_2(M)$, thus 
$E_{\omega_k}(\tilde{u}_k) = \int_{\dot{\Sigma}} u_k^*\omega_k =
\langle [\omega_k], [\bar{u}_k] \rangle = 0$.  This is a contradiction,
since $\tilde{u}_k$ is assumed to be nicely embedded, and thus nontrivial.
\end{proof}

\begin{lemma}
\label{lemma:avoidOrbits}
For each $k$, $u_k(\dot{\Sigma}) \subset M$ is disjoint from each of the
orbits $\gamma_{z,k} \subset M$ for $z \in \Gamma$.
\end{lemma}
\begin{proof}
Since $u_k$ is embedded, it follows from the nonlinear Cauchy-Riemann
equation that it is also transverse to $X_k$, thus any intersection
with $\gamma_{z,k}$ is transverse and implies transverse intersections of
$u_k$ with its image in a neighborhood of $z$.
\end{proof}

\begin{lemma}
\label{lemma:distinct}
For every $z \in \Gamma^\pm$, the extremal negative/positive eigenfunctions
of $\gamma_z$ with respect to $c_z$ are simply covered.  Moreover if $z$ 
and $\zeta$ are distinct
punctures with the same sign and $\gamma_z$ and $\gamma_\zeta$ cover the
same simply covered orbit, then $\gamma_z = \gamma_\zeta$.
\end{lemma}
\begin{proof}
Since $\asympdef(\tilde{u}_k ; \mathbf{c}) \le c_N(\tilde{u}_k ; \mathbf{c}) 
= 0$, $\tilde{u}_k$ is controlled by extremal eigenfunctions with respect to
$\mathbf{c}$ at each puncture, and these
must then be simply covered by Prop.~\ref{prop:uEmbedded}.  Similarly
Prop.~\ref{prop:uvDisjoint} implies that distinct positive/negative ends of 
$\tilde{u}_k$ approaching covers of the same orbit must approach with the
same covering number.  Both statements hold also in the limit due to the
nondegeneracy of the orbits $\gamma_z$.
\end{proof}

\begin{lemma}
\label{lemma:oppositeEnds}
Suppose $z_+ \in \Gamma^+$ and $z_- \in \Gamma^-$, such that $\gamma_{z_+}$ and
$\gamma_{z_-}$ cover the same simply covered orbit.  Then
$\gamma_{z_+} = \gamma_{z_-}$ and it is either a simply covered
even orbit or a doubly covered bad orbit with simply covered extremal
eigenfunctions.  Moreover, we can reset $c_{z_+} = c_{z_-} = 0$ without
changing $\ind(\tilde{u} ; \mathbf{c})$ or $c_N(\tilde{u} ; \mathbf{c})$.
\end{lemma}
\begin{proof}
For $\gamma_{z_\pm,k}$, the first part of the statement follows from 
Prop.~\ref{prop:oppositeSigns} since $u_k$ is embedded, and the second part
results from the fact that $\gamma_{z_\pm,k}$ is therefore even and
$\tilde{u}_k$ is controlled by extremal eigenfunctions (in the unconstrained
sense) at both of these punctures, so $\mathbf{A}_{\gamma_{z_\pm}}$ can
have no eigenvalues between $c_{z_\pm}$ and~$0$.  The same result is true for
$\gamma_{z_\pm}$ due to nondegeneracy.
\end{proof}

\begin{lemma}
\label{lemma:nonIntersecting}
For every connected component $\tilde{v}_i = (b_i,v_i) : \dot{S}_i 
\to \RR\times M$ of $\tilde{u}$, either $\tilde{v}_i$ is a trivial curve
or $v_i : \dot{S}_i \to M$ is injective.  Moreover for any two
such components $\tilde{v}_1$ and $\tilde{v}_2$ that are not trivial,
$v_1(\dot{S}_1)$ and $v_2(\dot{S}_2)$ are either disjoint or
identical, the latter if and only if $\tilde{v}_1$ can be obtained from
$\tilde{v}_2$ by $\RR$--translation (up to parametrization).
\end{lemma}
\begin{proof}
Suppose $\tilde{v}_i = (b_i,v_i) : \dot{S}_i \to \RR\times M$ is a
nontrivial connected component of $\tilde{u}$, so 
$E_\omega(\tilde{v}_i) > 0$ and consequently the section
$$
\pi Tv_i : \dot{S}_i \to \Hom_\CC(T\dot{S}_i,v_i^*\xi)
$$
has only finitely many zeros, all positive.  We claim first that $\tilde{v}_i$
is somewhere injective.  If not, then there exists a somewhere injective
finite energy surface $\tilde{w}_i = (\beta_i,w_i) : \dot{S}'_i \to
\RR\times M$ and a holomorphic branched cover $\varphi_i : \dot{S}_i \to
\dot{S}'_i$ of degree $k \ge 2$ such that $\tilde{v}_i =
\tilde{w}_i \circ \varphi_i$.  We can therefore find an embedded loop
$\alpha' : S^1 \to \dot{S}'_i$ which does not lift to $\dot{S}_i$, and
by small perturbations of $\alpha'$, we may assume it misses all punctures
and zeros of $\pi Tw_i$.  Now choose an embedded loop 
$\alpha : S^1 \to \dot{S}_i$ which projects to an $n$--fold cover of $\alpha'$
for some $n \ge 2$, and denote $C = \alpha(S^1) \subset \dot{S}_i$, 
$C' = \alpha'(S^1) \subset \dot{S}'_i$.  Choose also an open neighborhood
$\uU'$ of $C'$ and a corresponding neighborhood $\uU$ of $C$ such that
$\varphi_i(\uU) = \uU'$.  The restriction $\varphi_i|_{\uU} : \uU \to \uU'$
is an $n$--fold covering map, and we may assume without loss of generality
that $v_i|_{\uU} : \uU \to M$ and $w_i|_{\uU'} : \uU' \to M$ are both
transverse to~$X$.  From this we can derive a contradiction.
Indeed, any map $v' : \uU \to M$ that's $C^\infty$--close to $v_i|_{\uU}$
can be written on some neighborhood of $C$ as
$$
v'(z) = \varphi_X^{f(z)}(v_i(z))
$$
where $\varphi_X^t$ denotes the flow of $X$ and $f$ is a smooth real
valued function defined on some neighborhood of~$C$.  Choosing any
nontrivial deck transformation $g : \uU \to \uU$ for the covering map
$\varphi_i|_{\uU}$, there is necessarily a point $z \in C$ at which
$f(z) = f \circ g(z)$, and thus $v'(z) = v'(g(z))$.  By
Prop.~\ref{prop:convergence}, this is true
in particular for a suitable restriction of $u_k$ for~$k$ sufficiently
large, contradicting the assumption that $u_k$ is embedded.  We conclude
that $\tilde{v}_i$ is somewhere injective.

Now denote $\RR$--translations of finite energy surfaces $\tilde{u} = (a,u)$ by
$\tilde{u}^c := (a+c,u)$ for $c \in \RR$.
Suppose that $\tilde{v}_1$ and $\tilde{v}_2$ are two nontrivial
components and $v_1(z_1) = v_2(z_2)$.  This gives an intersection
$\tilde{v}_1(z_1) = \tilde{v}_2^c(z_2)$ for some $c \in \RR$.  If the
intersection is isolated then it is positive, and yields an isolated
intersection of $\tilde{u}_k$ and $\tilde{u}_k^{c'}$ for some $c' \in \RR$,
again contradicting the fact that $u_k$ is embedded.  The alternative,
since $\tilde{v}_1$ and $\tilde{v}_2$ are both somewhere injective, is
that $\tilde{v}_1$ and $\tilde{v}_2^c$ are identical up to parametrization.
The same argument applies to intersections of $v_1$ with itself:
since $\tilde{v}_1$ is somewhere injective,
the intersection of $\tilde{v}_1$ with $\tilde{v}_1^c$ is then necessarily
isolated, otherwise $\tilde{v}_1$ and $\tilde{v}_1^c$ are identical up to
parametrization; this is impossible in light of the asymptotic behavior
described in Prop.~\ref{prop:asymptotics}.
\end{proof}

We shall call a nonconstant trivial subbuilding of $\tilde{u}$
\emph{maximal} if every component attached to it by a breaking
orbit is nontrivial.  Given such a subbuilding 
$\tilde{u}\trivial : \mathbf{S}\trivial \to \RR\times M$, we introduce
the following notation: write the punctures of $\tilde{u}\trivial$ as
$$
\widehat{\Gamma}^\pm = \widehat{\Gamma}^\pm_{\oC} \cup \widehat{\Gamma}^\pm_{\oE},
$$
where $\widehat{\Gamma}_{\oE} := \widehat{\Gamma} \cap \Gamma$ and
$\widehat{\Gamma}_{\oC}$ consists of all punctures of $\tilde{u}\trivial$ 
that arise from \emph{breaking} punctures of $\tilde{u}$.
Assume $\#\widehat{\Gamma}^+_{\oC} = p$,
$\#\widehat{\Gamma}^-_{\oC} = q$, $\#\widehat{\Gamma}^+_{\oE} = r$ and
$\#\widehat{\Gamma}^-_{\oE} = s$; we have necessarily 
$\#\widehat{\Gamma}^+ = p + r > 0$,
$\#\widehat{\Gamma}^- = q + s > 0$ and since $\tilde{u}$ is connected and 
has nontrivial components, $\#\widehat{\Gamma}_{\oC} = p + q > 0$.  
Every asymptotic orbit of $\tilde{u}\trivial$ covers the same simply covered
orbit $\gamma$, so denote the orbit at $z \in \widehat{\Gamma}$ by
$$
\gamma_z = \gamma^{m_z}
$$
for some multiplicity $m_z \in \NN$.
Each $z \in \widehat{\Gamma}^\pm_{\oC}$ belongs to a breaking pair
$\{ z,\hat{z} \} \subset \Delta_{\oC}$ of $\tilde{u}$, and the component
$\tilde{v}_z = (b_z,v_z)$ of $\tilde{u}$ containing $\hat{z}$ is
necessarily nontrivial, and negatively/positively asymptotic to
$\gamma^{m_z}$ at~$\hat{z}$.  We know now from 
Lemma~\ref{lemma:nonIntersecting} that each of the maps $v_z$ is injective
(thus embedded near the punctures), and any two of them are either disjoint
or identical.  Then by the
intersection theoretic results of \S\ref{sec:intersection},
all the $m_z$ for $z \in \widehat{\Gamma}_{\oC}$ equal a fixed number
$m_{\oC} \in \NN$,
and the asymptotic approach of each $\tilde{v}_z$ to $\gamma_{m_{\oC}}$ is
controlled by eigenfunctions $e_z$ with the same winding
$\wind^\Phi(e_z) := w_{\oC} \in \ZZ$.  Likewise for $z \in \widehat{\Gamma}^\pm_{\oE}$,
Lemmas~\ref{lemma:distinct} and~\ref{lemma:oppositeEnds} imply that all 
$m_z$ equal a fixed multiplicity $m_{\oE} \in \NN$, and there is a fixed extremal 
winding number $w_{\oE} := \alpha^\Phi_\mp(\gamma^{m_{\oE}} ; c_z)$.  Note that
if both $\widehat{\Gamma}^+_{\oE}$ and
$\widehat{\Gamma}^-_{\oE}$ are nonempty, then
$\alpha^\Phi_+(\gamma^{m_{\oE}} ; c_z) = \alpha^\Phi_+(\gamma^{m_{\oE}}) =
\alpha^\Phi_-(\gamma^{m_{\oE}}) = \alpha^\Phi_-(\gamma^{m_{\oE}} ; c_z)$; 
this follows from Lemma~\ref{lemma:oppositeEnds}.

\begin{lemma}
\label{lemma:subbuildings}
For the maximal trivial subbuilding $\tilde{u}\trivial$ described above,
$m_{\oC} = m_{\oE}$ and $w_{\oC} = w_{\oE}$.
\end{lemma}
\begin{proof}
There's nothing to prove if $\widehat{\Gamma}_{\oE} = \emptyset$, so assume
$r + s > 0$.  Define the compact subset $\overline{\Sigma}\trivial = 
\psi(\overline{\mathbf{S}}\trivial)
\subset \overline{\Sigma}$ and recall that for sufficiently large $k$,
$\bar{u}_k \circ \psi|_{\overline{\mathbf{S}}\trivial}$ 
is $C^0$--close to $\bar{u}\trivial$.
Let $\gamma_k$ be the unique simply covered orbit of $X_k$ for sufficiently
large~$k$ such that $\gamma_k \to \gamma$.  Then some cover of $\gamma_k$
is an asymptotic orbit of $\tilde{u}_k$, thus by Lemma~\ref{lemma:avoidOrbits},
we can assume $\bar{u}_k(\overline{\Sigma}\trivial)$ 
lies in a fixed tubular neighborhood $N_\gamma$ of $\gamma$
but without intersecting $\gamma_k$.  We can also arrange that $\bar{u}_k$ 
have the following behavior at each component of $\p\overline{\Sigma}\trivial$:
\begin{itemize}
\item
for $z \in \widehat{\Gamma}^\pm_{\oC}$, 
$\wind^\Phi(\bar{u}_k(\psi(\delta_z))) = \pm w_{\oC}$,
\item
for $z \in \widehat{\Gamma}^\pm_{\oE}$, 
$\wind^\Phi(\bar{u}_k(\psi(\delta_z))) =
\pm w_{\oE}$.
\end{itemize}
The crucial observation is now that $\bar{u}_k(\overline{\Sigma}\trivial)$ 
realizes a homology
in $H_2(N_\gamma \setminus \gamma_k) \cong H_2(T^2)$.  From this we obtain
the relations
\begin{equation*}
\begin{split}
p m_{\oC} + r m_{\oE} &= q m_{\oC} + s m_{\oE}, \\
p w_{\oC} + r w_{\oE} &= q w_{\oC} + s w_{\oE},
\end{split}
\end{equation*}
and consequently
$$
\begin{pmatrix}
m_{\oC} & m_{\oE} \\
w_{\oC} & w_{\oE}
\end{pmatrix}
\begin{pmatrix}
p - q \\
r - s
\end{pmatrix} =
\begin{pmatrix}
0 \\
0
\end{pmatrix}.
$$
If $p=q$ and $r=s$, then both of these are nonzero and 
Prop.~\ref{prop:oppositeSigns} implies that either $\gamma$ is even and
$m_{\oC} = m_{\oE} = 1$ or $\gamma$ is odd hyperbolic and $m_{\oC} = m_{\oE} = 2$,
with $w_{\oC} = w_{\oE} = \alpha^\Phi_+(\gamma^{m_{\oC}}) = \alpha^\Phi_-(\gamma^{m_{\oC}})$ 
in either case. 
Otherwise the determinant of the matrix above must vanish, so
$m_{\oC} w_{\oE} = m_{\oE} w_{\oC}$.  However, by Prop.~\ref{prop:uEmbedded},
$w_{\oC}$ and $w_{\oE}$ are winding numbers of simply covered eigenfunctions 
for $\gamma^{m_{\oC}}$ and $\gamma^{m_{\oE}}$ respectively, thus
Prop.~\ref{prop:relPrime} implies that $m_{\oC}$ and $w_{\oC}$ are relatively
prime, as are $m_{\oE}$ and $w_{\oE}$.  This implies $m_{\oC} = m_{\oE}$ and
$w_{\oC} = w_{\oE}$
\end{proof}

\begin{cor}
\label{cor:subbuildings}
If $\tilde{u}\trivial : \mathbf{S}\trivial 
\to \RR\times M$ is the maximal trivial 
subbuilding above with induced asymptotic constraints $\hat{\mathbf{c}} = 
\{ \hat{c}_z \}_{z \in \widehat{\Gamma}}$, then
$$
c_N(\tilde{u}\trivial ; \hat{\mathbf{c}}) +
\sum_{z \in \widehat{\Gamma}_{\oC}}\left[ p(\gamma_z) +
\asympdef^{\hat{z}}(\tilde{v}_z) \right] = 
- \chi(\overline{\mathbf{S}}\trivial).
$$
In particular this sum is nonnegative, and is zero if and only if
$\tilde{u}\trivial$ is cylindrical.
\end{cor}
\begin{proof}
By the lemma we have $\#\widehat{\Gamma}^+ = \#\widehat{\Gamma}^-$ and 
can write $m := m_{\oE} = m_{\oC}$ and $w := w_{\oE} = w_{\oC} = 
\alpha^\Phi_\mp(\gamma^m ; c_z)$ for each $z \in \widehat{\Gamma}^\pm_{\oE}$.  
Then, noting that $c_1^\Phi((\bar{u}\trivial)^*\xi) = 0$,
\begin{equation*}
\begin{split}
c_N(\tilde{u}\trivial ; \hat{\mathbf{c}}) &+
\sum_{z \in \widehat{\Gamma}_{\oC}} \left[ p(\gamma_z) +
\asympdef^{\hat{z}}(\tilde{v}_z) \right] \\
&= -\chi(\overline{\mathbf{S}}\trivial) + 
\sum_{z \in \widehat{\Gamma}^+} \alpha^\Phi_-(\gamma^m ; \hat{c}_z)
 - \sum_{z \in \widehat{\Gamma}^-} \alpha^\Phi_+(\gamma^m ; \hat{c}_z) \\
&\qquad + \sum_{z \in \widehat{\Gamma}_{\oC}} \left[ 
\alpha^\Phi_+(\gamma^m) - \alpha^\Phi_-(\gamma^m) \right] \\
&\qquad + \sum_{z \in \widehat{\Gamma}^+_{\oC}} 
\left[ w - \alpha^\Phi_+(\gamma^m) \right] +
\sum_{z \in \widehat{\Gamma}^-_{\oC}} \left[ \alpha^\Phi_-(\gamma^m) - w \right] \\
&= -\chi(\overline{\mathbf{S}}\trivial) +
\sum_{z \in \widehat{\Gamma}^+} w -
\sum_{z \in \widehat{\Gamma}^-} w 
= -\chi(\overline{\mathbf{S}}\trivial).
\end{split}
\end{equation*}
\end{proof}

We shall handle constant components of $\tilde{u}$ similarly.
Call a connected subbuilding $\tilde{u}\const : \mathbf{S}\const 
\to \RR\times M$ 
of $\tilde{u}$ a \emph{constant subbuilding} if every connected component of 
$\tilde{u}\const$ is constant.  Further, call it a \emph{maximal} constant 
subbuilding if every constant component of $\tilde{u}$ that is attached 
by a node to some component of $\tilde{u}\const$ 
is also in $\tilde{u}\const$.  Note
that constant subbuildings cannot have punctures, thus 
$\overline{\mathbf{S}}\const$ is closed.  
Denote by $\widehat{\Delta}_{\oN} \subset
S\const$ the set of points $z \in S\const$ that belong to nodal pairs
$\{z,z'\} \subset \Delta_{\oN}$ of $\mathbf{S}$ such that $\tilde{u}$ is not
constant near $z'$; this set is necessarily nonempty since
$\tilde{u}$ is connected.  Then the stability condition on $\tilde{u}$ implies
$$
\chi(\overline{\mathbf{S}}\const \setminus \widehat{\Delta}_{\oN}) < 0.
$$
Thus $c_N(\tilde{u}\const) + 2 \#\widehat{\Delta}_{\oN} = 
-\chi(\overline{\mathbf{S}}\const) + \#\widehat{\Delta}_{\oN} + 
\#\widehat{\Delta}_{\oN}
= -\chi(\overline{\mathbf{S}}\const \setminus \widehat{\Delta}_{\oN}) +
\#\widehat{\Delta}_{\oN} > \#\widehat{\Delta}_{\oN} > 0$.  We've proved:

\begin{lemma}
\label{lemma:constant}
For any maximal constant subbuilding $\tilde{u}\const$ of $\tilde{u}$ with
nodes $\widehat{\Delta}_{\oN}$ connecting it to nonconstant components of
$\tilde{u}$,
$$
c_N(\tilde{u}\const) + 2 \#\widehat{\Delta}_{\oN} > 0.
$$
\end{lemma}

All the ingredients are now in place.

\begin{proof}[Proof of Theorem~\ref{thm:mainTechnical}]
By the above results, $\tilde{u}$ consists of the following pieces:
\begin{enumerate}
\item
Maximal constant subbuildings $\tilde{u}\const$ such that
$c_N(\tilde{u}\const) + 2 \#\widehat{\Delta}_{\oN} > 0$.
\item
Maximal trivial subbuildings $\tilde{u}\trivial$ with induced asymptotic
constraints $\mathbf{c}\trivial$, for which the sum of
$c_N(\tilde{u}\trivial ; \mathbf{c}\trivial) + 
\sum_{z \in \widehat{\Gamma}_{\oC}} p(\gamma_z)$ with the asymptotic
defects of all neighboring nontrivial ends is nonnegative, and zero if and
only if $\tilde{u}\trivial$ is cylindrical.
\item
Nontrivial connected components $\tilde{v} = (b,v)$ with $v$ injective.
\end{enumerate}
Note that each nontrivial component $\tilde{v}$ is compatible with induced
asymptotic constraints $\hat{\mathbf{c}}$ and satisfies
$c_N(\tilde{v} ; \hat{\mathbf{c}}) - \asympdef(\tilde{v} ; \hat{\mathbf{c}})
= \windpi(\tilde{v}) \ge 0$.

Since $c_N(\tilde{u} ; \mathbf{c}) = 0$, we conclude from
Prop.~\ref{prop:cNsum} that $\tilde{u}$ contains no constant
subbuildings or nodes, every trivial subbuilding is cylindrical and 
every nontrivial component $\tilde{v}$ has $\windpi(\tilde{v}) = 0$.
Such components $\tilde{v}$ are therefore nicely embedded.
A slightly stronger statement results from the observation
that the core $\tilde{u}\core$ is also compatible with $\mathbf{c}$
and only contains nicely embedded components.  Thus each of these components
$\tilde{v}$ satisfies $c_N(\tilde{v} ; \hat{\mathbf{c}}) = 0$, where
$\hat{\mathbf{c}}$ are now 
the constraints induced on $\tilde{v}$ as a subbuilding of~$\tilde{u}\core$.
\end{proof}

\begin{proof}[Proof of Theorem~\ref{thm:stableTechnical}]
For a given set of asymptotic constraints $\mathbf{c}$, the $\RR$--invariance
of $\tilde{J}$ together with a standard
transversality argument (cf.~\cite{Wendl:BP1}) imply that for generic 
$\omega$--compatible choices of $J$, all nontrivial
somewhere injective finite energy surfaces $\tilde{w}$ compatible with 
$\mathbf{c}$ satisfy $\ind(\tilde{w} ; \mathbf{c}) \ge 1$.  Moreover,
$c_N(u ; \mathbf{c}) = 0$ implies
$\ind(u ; \mathbf{c}) \le 2$ due to \eqref{eqn:cNindex}, thus
this index can only be~$1$ or~$2$.  Likewise each
connected component of $\tilde{u}\core$ has constrained
index at least~$1$, and by Prop.~\ref{prop:indexAdditivity}, these add
up to $\ind(\tilde{u} ; \mathbf{c})$.  We conclude there is exactly
one component if $\ind(\tilde{u} ; \mathbf{c}) = 1$, and at most two if
$\ind(\tilde{u} ; \mathbf{c}) = 2$.  In the latter case, both
nontrivial components $\tilde{v}_i$ have 
$\ind(\tilde{v}_i ; \mathbf{c}_i) = 1$,
so by \eqref{eqn:cNindex}, each has a unique puncture whose constrained
parity is even: this is therefore the unique breaking puncture.
Since the ends of $\tilde{v}_1$ and $\tilde{v}_2$ approaching this
breaking orbit have opposite signs, $\tilde{v}_1$ and $\tilde{v}_2$
cannot be the same up to $\RR$--translation, thus their projections to~$M$
are disjoint.
\end{proof}

\end{document}